\documentclass{article}

\usepackage[preprint]{neurips_2026}

\usepackage{amsmath,amsfonts,bm}

\def\Figref#1{Figure~\ref{#1}}

\def\twoFigref#1#2{Figures \ref{#1} and \ref{#2}}

\def\Secref#1{Section~\ref{#1}}
\def\eqref#1{equation~\ref{#1}}
\def\plaineqref#1{\ref{#1}}
\def\1{\bm{1}}

\def\rvn{{\mathbf{n}}}

\def\mA{{\bm{A}}}

\def\mF{{\bm{F}}}

\def\mH{{\bm{H}}}

\def\mM{{\bm{M}}}

\def\mU{{\bm{U}}}
\def\mV{{\bm{V}}}

\DeclareMathAlphabet{\mathsfit}{\encodingdefault}{\sfdefault}{m}{sl}
\SetMathAlphabet{\mathsfit}{bold}{\encodingdefault}{\sfdefault}{bx}{n}

\def\sC{{\mathbb{C}}}

\def\sN{{\mathbb{N}}}

\def\sR{{\mathbb{R}}}

\DeclareMathOperator*{\argmin}{arg\,min}
\DeclareMathOperator{\prox}{prox}

\usepackage{xspace}

\newcommand{\method}{i-DEQ\xspace}
\newcommand{\methodP}{i-DEQ-P\xspace}

\newcommand{\Dsf}{\mathsf{D}}

\newtheorem{theorem}{Theorem}[section]

\newtheorem{proposition}{Proposition}[section]

\newtheorem{remark}{Remark}
\newtheorem{assumption}[theorem]{Assumption}
\usepackage[utf8]{inputenc}
\usepackage[T1]{fontenc}

\usepackage{hyperref}

\usepackage{booktabs}
\usepackage{amsfonts}
\usepackage{mathtools}
\usepackage{nicefrac}

\usepackage{microtype}
\usepackage{xcolor}
\usepackage{xspace}

\usepackage{graphicx}

\usepackage{tikz}
\usetikzlibrary{spy,calc}

\usepackage{array}
\usepackage{caption}
\usepackage[normalem]{ulem}

\usepackage{algorithm}
\usepackage{algpseudocode}

\usepackage{multirow}
\usepackage{enumitem}
\usepackage{wrapfig}

\usepackage{placeins}
\graphicspath{{pics/}}

\title{\method: A stable inertial deep equilibrium model for image restoration}

\author{%
  Antonin Clerc \thanks{corresponding author: \texttt{antonin.clerc@math.u-bordeaux.fr}} \\
  Univ. Bordeaux, CNRS, Inria, Bordeaux INP, IMB, UMR 5251\\
  F-33400 Talence, France \\
  \And
  Marien Renaud \\
  Univ. Bordeaux, CNRS, Bordeaux INP, IMB, UMR 5251 \\
  F-33400 Talence, France \\
  \AND
  Baudouin Denis de Senneville \\
  Univ. Bordeaux, CNRS, Inria, Bordeaux INP, IMB, UMR 5251\\
  F-33400 Talence, France \\
  \And
  Nicolas Papadakis \\
  Univ. Bordeaux, CNRS, Inria, Bordeaux INP, IMB, UMR 5251\\
  F-33400 Talence, France \\
}

\definecolor{forbg}{RGB}{245, 235, 200}
\definecolor{backbg}{RGB}{220, 240, 220}

\begin{document}

\maketitle

\begin{abstract}
Deep Equilibrium Models (DEQs) are an established framework for image restoration that learn a problem-adapted regularization by solving a fixed-point (i.e. equilibrium) problem.
While flexible and expressive, DEQs are often hindered by high computational cost and training instability.
We propose an inertial DEQ (\method) that learns an explicit nonconvex regularization within the DEQ formulation.
By using momentum within the fixed-point iterations, \method has convergence guarantees and accelerated rates. 
Moreover, we observe that \method is significantly more stable during the training and robust to rough initialization than DEQs.
Numerical experiments on various linear and nonlinear 
inverse problems demonstrate that \method achieves reconstruction quality comparable to state-of-the-art methods, while reducing DEQ's inference time by a factor of 2 to 5.
\end{abstract}

\section{Introduction}
Inverse problems in imaging aim to recover an unknown image $x \in \sR^d$ from noisy and degraded observations $y \in \sR^m$, modeled as $y = \mathcal{A}(x) + \rvn$, where $\mathcal{A}:\sR^d\to\sR^m$ is a forward operator and $\rvn$ is typically an additive white Gaussian noise (AWGN). Arising in applications from computer vision to biomedical imaging~\citep{bertero1998introduction,tarantola2005inverse},
these problems are ill-posed~\citep{hadamard1902problemes}, meaning that the observations do not uniquely determine the unknown image. To address this issue, model-based approaches incorporate prior knowledge through regularization~\citep{engl1996regularization}. A standard formulation is Maximum a Posteriori (MAP) estimation, which seeks
\begin{equation} \label{eq:inv_prob}
    \hat{x} \in \arg\min_{x \in \sR^d} F(x) := f(x,y) + \lambda g(x),
\end{equation}
where $f$ measures the fidelity to the data $y$ and $g$ encodes prior information on the image to recover.
Many explicit priors have been proposed, including sparsity in transform domains~\citep{donoho2006compressed,mallat1999wavelet}, total variation~\citep{rudin1992nonlinear}, low-rank models~\citep{candes2011robust}, and non-local self-similarity~\citep{buades2005nonlocal}.

Leveraging the expressiveness of deep neural networks, Model-Based Deep Learning (MBDL) has emerged, in which priors are no longer hand-crafted but learned from data~\citep{monga2021algorithm}.
Several classes of MBDL methods have been proposed, see a recent review in~\citep{hertrich_learning_2025}. 
In particular, Plug-and-Play (PnP)~\citep{venkatakrishnan2013plug} approaches integrate a denoiser within iterative optimization algorithms. 
However, since the regularization imposed by the denoiser is implicit, establishing convergence guarantees becomes challenging and requires strong assumptions on the underlying operators \citep{combettes2005signal}.
Therefore, without adequate safeguards, the execution of a PnP algorithm can result in divergence throughout the restoration procedure~\citep{nair2024averaged,Reehorst2019,sommerhoff2019energy,terris2024equivariant,zhang2017beyond}.
In contrast, PnP methods relying on explicit regularization define a well-defined optimization problem~\citep{HuraultGSD,romano2017red}, thereby facilitating theoretical analysis and convergence guarantees. 
Nevertheless, MBDL approaches typically require a large number of iterations to converge, and their performance depends on careful parameter tuning.

Rather than explicitly designing a regularization functional, End-to-End (E2E) deep learning learns a direct mapping from measurements to images, in a purely data-driven manner~\citep{jin2017deep}. These methods achieve state-of-the-art performance in many imaging tasks \citep{RAM}. While highly effective, they can disregard the underlying acquisition model and do not explicitly enforce data-consistency, which may limit robustness and generalization~\citep{Antun2020}.
Unrolled approaches map a finite number of iterations of an optimization algorithm onto a neural network trained E2E for a specific forward model~\citep{chen2018theoretical,gregor2010learning}. This leads to improved reconstruction performance by incorporating domain knowledge into the architecture \citep{monga2021algorithm}. However, these methods suffer from high memory consumption during training, as backpropagation requires storing intermediate states across the unrolled iterations, leading to a cost that increases with the number of layers. As a result, this memory cost limits the depth to a small number of iterations, typically fewer than 10, and restricts the scalability.
In addition, they typically lack an explicit optimization objective, making theoretical analysis more involved.

Deep Equilibrium (DEQ)~\citep{bai_multiscale_2020,Gilton2021DEQImaging} models bridge the gap between E2E approaches and MBDL ones. They inherit their optimization structure from MBDL, but unlike 
PnP methods, which rely on off-the-shelf denoisers, DEQs learn the underlying operator directly from the data of the tackled problem. This yields problem-adapted models that can match or even surpass E2E approaches \citep{ghaoui_implicit_2020,ling2023deep,sun2023understanding,tsuchida2022deep,wu2023separation}, while still preserving a variational interpretation, in contrast to standard E2E~\citep{DEQ_Poisson,Zou2023ELDER}. 
Moreover, DEQs are more memory-efficient to train than unrolled methods, as they rely on implicit differentiation and avoid backpropagation through all intermediate iterations~\citep{bai_multiscale_2020,JFB}. However, computing the solution requires solving a fixed-point problem, which remains computationally demanding despite the use of acceleration techniques or black-box fixed-point solvers~\citep{bai2021deq,lin_consistency_2026, winston_monotone_2021}. In addition, training DEQ models is often prone to instability, and designing effective strategies to stabilize the training process remains challenging~\citep{amos2021optnet,bai2021stabilizing,blondel2022efficient,geng2022training, mccallum2025reversible}.

\paragraph{Contributions}
In this work, we introduce an accelerated tuning-free DEQ-based method that stabilizes the training by incorporating restarted inertia in the optimization iterations.

\begin{figure}[htbp] 
    \begin{tabular}{cc}
    \includegraphics[width=0.48\linewidth]{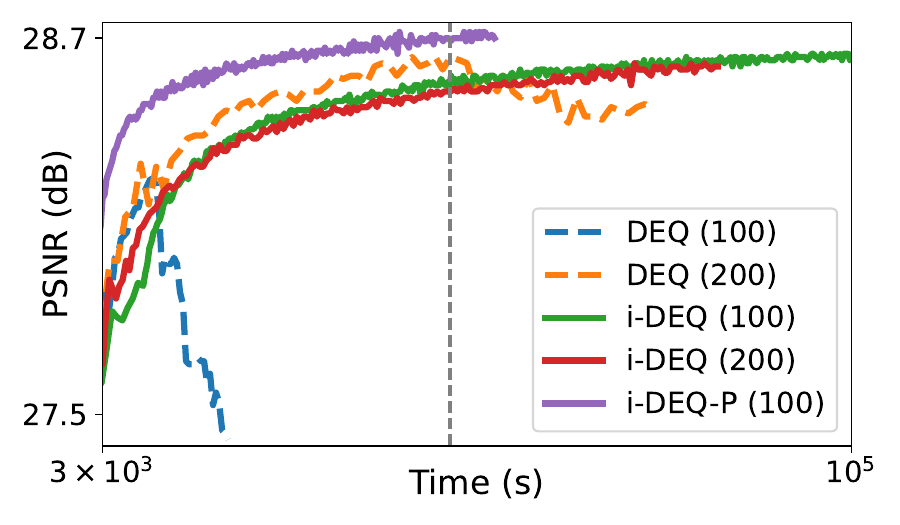} &
    \includegraphics[width=0.48\linewidth]{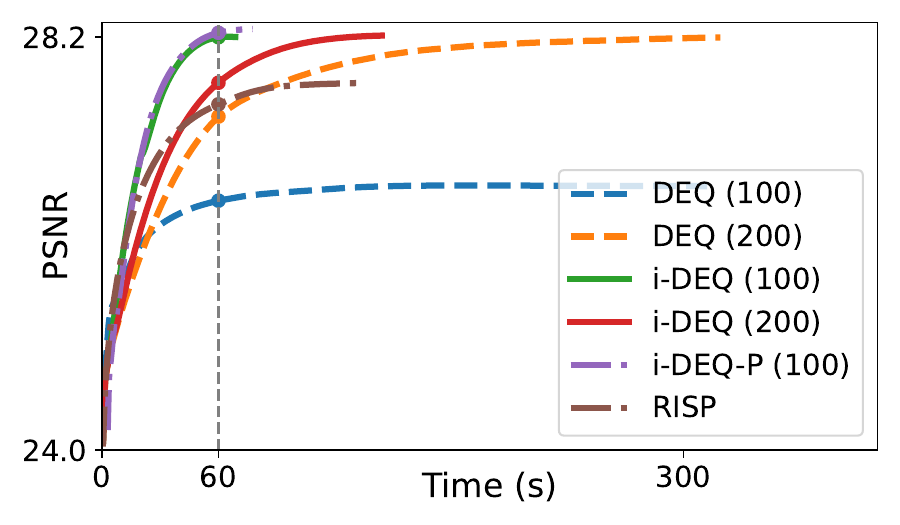}\vspace{-.2cm}\\
    (a) Training&(b) Inference
    \end{tabular}
    \caption{Training on 100 images and inference on a test set of 20 images for DEQ-based methods on a subset of FastMRI \citep{knoll2020fastmri}. 
    \textbf{(a)} impact of the number of fixed-point iterations in the forward pass during training. Standard DEQ becomes unstable for both 100 and 200 iterations, whereas \method remains stable and shows a steady PSNR improvement. \method-P (100) exhibits faster improvement than \method and DEQ while being more stable than DEQ.  
    \textbf{(b)} PSNR vs inference time for reaching the stopping criterion. \method, trained with 100 and 200 iterations, is compared to a classical DEQ with matching iteration budgets and to RISP~\citep{RISP}; it achieves similar PSNR  2--5$\times$ faster.}
    \label{fig:metric_vs_time}
    \vspace{-0.3cm}
\end{figure}

\begin{itemize}
    \item By relying on an explicit regularization formulation derived from PnP, the method enjoys well-definedness guarantees, including the existence of stationary points and the asymptotic convergence of the iterative scheme toward critical points of $F$ under standard assumptions (see Proposition~3.1).
    
    \item We leverage inertial dynamics within the optimization iterations to accelerate convergence toward critical points during training, which empirically leads to improved training stability compared to classical DEQ formulations (see \Figref{fig:metric_vs_time} left).
    
    \item The use of an inertial acceleration scheme also significantly reduces the number of iterations required to reach equilibrium at inference, enabling faster reconstructions without loss of performance (see \Figref{fig:metric_vs_time} right).\vspace*{-.3cm}
\end{itemize}

\paragraph{Outline of the paper}
In \Secref{sec:related_work}, we introduce some useful tool and survey related approaches in data-driven reconstruction methods.
In \Secref{sec:method}, we present the proposed inertial Deep Equilibrium (\method) and provide its theoretical accelerate rate.
In \Secref{sec:experiments}, we validate our method on inverse problems including compressed sensing MRI, inpainting, motion deblurring and non-linear Rician denoising. In particular, we show competitive performances with state-of-the-art approaches, while exhibiting improved stability during training. In addition, we demonstrate faster restoration time compared to PnP and standard DEQs.
Code, data, and network weights are provided to the reviewers and will be made publicly available upon publication. Finally, conclusions are drawn in Section 5.

\section{Relationship to Prior Work} \label{sec:related_work}

In this section, we review the main approaches related to our work. We first recall classical first-order optimization methods. We then discuss Plug-and-Play frameworks, which incorporate denoisers as implicit priors within iterative reconstruction schemes. Finally, we present learning-based reconstruction methods, including end-to-end and deep equilibrium models, which learn a reconstruction operator that is directly adapted to the inverse problem.

\subsection{Optimization}

Gradient descent (GD) is a standard iterative method for solving problem~\plaineqref{eq:inv_prob}~\citep{cauchy1847methode}. Given an initial estimate $x^0$ and a step size $\tau > 0$, the iterates are defined as $x^{(k+1)} = x^{(k)} - \tau \nabla F(x^{(k)})$. 

In the case of a composite objective $F = f + g$ with a smooth term $f$ and a possibly non-smooth term $g$, proximal gradient descent (PGD) can be used~\citep{combettes2011proximal,parikh2014proximal}. It alternates a gradient step on $f$ with a proximal step associated with $g$, where the proximal operator of $g : \sR^d \to \sR$ is defined as
$\mathrm{prox}_{\tau g}(x) 
= \argmin_{u \in \sR^d} 
\left\{
g(u) + \frac{1}{2\tau}\|u - x\|^2
\right\}$, 
which  leads to the  iterations  
$x^{k+1} = \operatorname{prox}_{\tau g}\big(x^k - \tau \nabla f(x^k,y)\big)$.

\subsection{Plug-and-Play and Regularization by Denoising Methods}
Plug-and-Play (PnP) and Regularization by Denoising (RED) are iterative algorithms for solving inverse problems of the form~\plaineqref{eq:inv_prob}, with a prior $g$ encoded by a denoiser. 
PnP was introduced as forward-backward schemes where an off-the-shelf denoiser is plugged in place of the proximal operator of the regularization term~\citep{venkatakrishnan2013plug}.
However, in general, the underlying objective function is not explicitly defined unless additional assumptions are imposed on the denoiser~\citep{combettes2020lipschitz,  hertrich2021convolutional, HuraultGSD,pesquet_learning_2021}. 
Alternatively, RED uses an off-the-shelf denoiser to compute the gradient of an explicit regularization~\citep{romano2017red}. More precisely, a denoiser $\Dsf_\sigma : \sR^d \rightarrow \sR^d$, with parameter $\sigma > 0$ controlling the denoising strength, is used to construct a gradient-like quantity through the residual $I_d - \Dsf_\sigma$. Restrictive assumptions such as symmetric Jacobian~\citep{romano2017red} are still required on the denoiser for the explicit regularizer to exist.

\paragraph{Gradient Step Denoising}
The authors of~\citep{HuraultGSD} ensures the convergence of RED by introducing an explicit regularization of the form
\begin{equation} \label{eq:GS}
\Dsf_\sigma = \mathrm{Id} - \nabla g_\sigma, 
\qquad 
g_\sigma(x) \coloneq \frac{1}{2}\|x - N_\sigma(x)\|^2.
\end{equation}

Here, $g_\sigma$ is a learnable potential parameterized by a neural network $N_\sigma$, allowing it to retain the expressive power of modern learned denoisers. By construction, $g_\sigma$ is continuously differentiable, so that $\Dsf_\sigma$ can be interpreted as the gradient of a scalar potential, ensuring consistency with a variational formulation.
The resulting RED update with step size $\tau>0$ is given by:
\begin{equation}\label{eq:RED_operator} 
\textbf{RED:} \quad x^{k+1} = T^{\mathrm{RED}}(x^k) 
\coloneqq x^k - \tau \left( \nabla f(x^k,y) + \lambda \big(x^k - \Dsf_\sigma(x^k)\big) \right).
\end{equation}

% where $\theta$ denotes the set of algorithmic parameters, i.e. $\theta = \{\tau, \lambda, \sigma\}$, and $\lambda$ controls the influence of the denoiser $D_\sigma$.

\paragraph{Accelerating RED}
The Restarted Inertia and Score-based Image Priors method (RISP)~\citep{RISP} proposes an acceleration of RED algorithms by combining inertial extrapolation with a restarting mechanism. Given an iterate $ x^k $, the extrapolated point is defined as
\begin{equation} \label{eq:Momentum}
z^k = T^\text{M}_\alpha(x^k, x^{k-1}) \coloneq x^k + (1-\alpha)(x^k - x^{k-1}),
\end{equation}
where $ \alpha \in (0,1] $ controls the strength of the inertial term.
This extrapolated point is taken as input of the RED operator \plaineqref{eq:RED_operator}, so that RISP corresponds to RED for $\alpha=1$.

Nevertheless, accumulated inertia may induce instability and undesirable overshooting, particularly in nonconvex settings~\citep{odono2015adaptive}. To mitigate this issue, RISP incorporates a restarting mechanism~\citep{li2023restarted} that resets the inertial contribution whenever the accumulated variation exceeds a prescribed threshold parameterized by the user-defined constant $B > 0$
\begin{equation} \label{eq:restart}
    k\sum_{t=0}^{k-1} \|x^{t+1} - x^t\|^2 > B^2.
    %\vspace{-0.3cm}
\end{equation}

\subsection{Problem adaptive reconstruction methods}

\paragraph{Deep unrolled networks}
Unrolled optimization methods build deep networks by truncating an iterative algorithm to a fixed number of steps, each mapped to a learnable block. They achieve competitive performances on imaging inverse problems~\citep{chun2020momentumnet,galinier2020hybrid, gharbi_unrolled_2024,mehranian2020model, monga_algorithm_2020,wu2019computationally}
.
Unrolling interprets $K$ iterations as a feedforward network trained by backpropagating through all the layers. While gradients are exact for the truncated procedure, this requires storing all intermediate states and operations, which limits $K$ (typically 5–10 in MRI~\citep{pramanik_adapting_2023}).
Since the architecture consists of a sequence of similar blocks, a key design choice is whether to use independent parameters at each iteration or to share them across iterations (\emph{weight sharing/tying}). In practice, parameters are often not shared~\citep{schlemper_deep_2017}, increasing flexibility but moving these models away from classical iterative schemes.

\paragraph{Deep Equilibrium Networks}

Deep Equilibrium (DEQ) models define the reconstructed image as the fixed-point of a learned nonlinear operator. This operator is parameterized by a neural network with \emph{shared weights} and applied iteratively until convergence, rather than being unrolled for a fixed number of iterations.
As a result, DEQs can be interpreted as implicit infinite-depth networks, combining the expressivity of E2E architectures, since they encompass them as a special case~\citep{bai_multiscale_2020,ghaoui_implicit_2020,tsuchida2022deep}, with the structured nature of MBDL approaches.
These models have been widely applied to inverse problems in imaging, including MRI~\citep{gan_selfsupervised, Gilton2021DEQImaging,mehta_equivariant_2025,yang_equilibrium_2023,Zou2023ELDER}, computed tomography~\citep{Liu2022ODEQ,chazottes_deep_2024}, video reconstruction~\citep{Zhao2022DEQVideo}, and nonlinear inverse problems such as Poisson denoising~\citep{DEQ_Poisson}.

For $T_\theta : \sR^d \rightarrow \sR^d$ a parameterized operator, a DEQ model defines its output as the equilibrium point
$
\bar x = T_\theta(\bar x).
$
In practice, this fixed-point is approximated using the iterative scheme
\begin{equation*}
    x^{k+1} = T_\theta(x^k), \quad k \geq 0,
\end{equation*}
until convergence, which requires a large number of iterations to compute the equilibrium point.
From a theoretical perspective, the well-posedness of DEQ models relies on the existence of a fixed-point, as well as the convergence of the  iterates $x^k$~\citep[Definition~1]{DEQ_Poisson} to a fixed-point $\bar x$, properties  related to classical conditions in model-based frameworks such as PnP and RED. Early DEQ reconstruction methods~\citep{Gilton2021DEQImaging} leverage implicit PnP  regularization, making convergence analysis difficult unless strong constraints are imposed on the denoiser, e.g.  with monotone operators~\citep{pramanik_improved_2021}.
The ELDER method~\citep{Zou2023ELDER} overcomes this limitation by introducing an explicit regularization based on the GS-Denoiser~\plaineqref{eq:GS}, leading to the minimization of an explicit energy functional with convergence guarantees.
However, it requires backtracking to select step sizes, which increases computational cost due to repeated evaluations of the energy and recomputation of fixed-point iterations.

DEQ training relies on implicit differentiation of the fixed-point equation, avoiding backpropagation through the iterative solver~\citep{bai_multiscale_2020}, most commonly via Jacobian-Free Backpropagation (JFB)~\citep{JFB}.
This yields constant memory complexity compared to unrolled networks. The  procedure consists of:
\begin{itemize}[leftmargin=12pt]
    \item[1.] \textbf{Forward pass:} approximate fixed-point using an iterative solver or root-finding algorithm.
    \item[2.]\textbf{Backward pass:} compute gradients via implicit differentiation of the equilibrium equation.
\end{itemize}
An epoch corresponds to one computationally demanding full forward fixed-point computation, followed by a backward pass based on implicit differentiation of a loss function $\mathcal{L}$, typically the MSE loss,  between the computed fixed-point $\bar x$ and the true signal $x^\star$ (see Appendix~\ref{app:training_DEQ_gen} for more details).
DEQ models are known to suffer from numerical instabilities during backpropagation, and prior works address this issue either through architectural constraints~\citep{amos2021optnet,bai2021stabilizing, mccallum2025reversible} or refined implicit differentiation schemes~\citep{blondel2022efficient, davy_restarted_2025,geng2022training}, which introduces additional computational overhead.
Thus, despite their strong theoretical appeal, DEQ models remain challenging to deploy in practice, motivating the development of more stable and computationally efficient equilibrium-based formulations that preserve the variational structure while reducing the cost of fixed-point iterations.

\section{Method} \label{sec:method}

We introduce \method, a novel DEQ framework for inverse problems that learns the potential function $g_\theta$ through energy minimization. 
Unlike standard DEQ approaches~\citep{Gilton2021DEQImaging} defined via implicit fixed-point iterations, our forward pass explicitly minimizes a learned energy functional. 
This variational perspective enables the use of inertially accelerated optimization with formal convergence guarantees, combining the interpretability of energy-based priors with stable, efficient execution. 
In what follows, we outline the \method architecture and establish its theoretical convergence.

\subsection{Energy-Based Forward Pass}

We consider inverse problems of the form~\plaineqref{eq:inv_prob} and aim to learn a problem-adapted regularization functional while preserving a variational formulation. The regularizer is parameterized as a scalar energy $g_\theta$ and built using the GS-Denoiser structure defined in~\plaineqref{eq:GS}.
The forward pass consists in finding the fixed-point of a learned inertial optimization operator. Starting from the RED formulation~\plaineqref{eq:RED_operator} and leveraging the explicit regularization induced by the GS-Denoiser as in~\citep{Zou2023ELDER}, together with the inertial update~\plaineqref{eq:Momentum}, we define the combined operator:
\begin{align}
T_{\theta,\lambda,\tau}^{\mathrm{RED}}(x^k)
&\coloneq x^k - \tau\Big(\nabla f(x^k,y) + \lambda \nabla g_\theta(x^k)\Big), \nonumber\\
T_\Theta^{\mathrm{\method}}(x^k, x^{k-1})
&\coloneq T_{\theta,\lambda,\tau}^{\mathrm{RED}}\circ T_\alpha^{\mathrm{M}}(x^k, x^{k-1}),\label{eq:ideq}
\end{align}

\begin{wrapfigure}{r}{0.58\textwidth} % r = droite, 0.48 = largeur (ajustable)
\vspace{-0.4cm} % Remonte l'algo pour l'aligner avec la première ligne du texte

\begin{minipage}{\linewidth}
\footnotesize
\setlength{\baselineskip}{10pt}

% Reconstitution propre de l'en-tête d'algorithme sans le flottant principal
\rule{\linewidth}{0.8pt}\vspace{-3pt}
\captionof{algorithm}{\method}
\label{algo:GM_RISP_DEQ}
\vspace{-6pt}\rule{\linewidth}{0.4pt}\vspace{2pt}

\setlength{\fboxsep}{1pt}
\begin{algorithmic}[1]
\Require Dataset $(Y,X^\star)$, $n_{ep}, K \in \sN$, $\lambda, \tau > 0$, $\alpha \in (0,1)$, $B > 0$

\For{$n = 0$ to $n_{ep}-1$}

\State \hspace{-\fboxsep}\colorbox{forbg}{\parbox{\dimexpr\linewidth-\algorithmicindent-1pt}{$\nabla_\Theta L\gets 0$\hfill\textbf{Forward pass} }}

\For{each $(y,x^\star) \in (Y,X^\star)$} 
    \State $k \gets 0$, \quad $x^{-1} = x^0 \in \sR^d$
    
    \While{$k \leq K$}
    
        \State \hspace{-\fboxsep}\colorbox{lightgray!50!white}{\parbox{\dimexpr\linewidth-27pt-\algorithmicindent-2\fboxsep}{
        $z^k = x^k + (1 - \alpha)(x^k - x^{k-1})$
        \hfill \textcolor{black}{\textit{Inertial step} }
        }}
        
        \State $x^{k+1} \gets z^k - \tau \left( \nabla f(z^k, y) + \lambda \nabla g_\theta(z^k) \right)$
        
        \State $k \gets k + 1$
        
        \If{$k \sum_{i=0}^{k-1} \|x^{i+1} - x^i\|^2 > B^2$}
            \State \hspace{-\fboxsep}\colorbox{lightgray!50!white}{\parbox{\dimexpr\linewidth-40pt-\algorithmicindent-2\fboxsep}{
            $x^{-1}, x^0 \gets x^k$, \quad $k \gets 0$
            \hfill \textcolor{black}{\textit{Restart} }
            }}
        \EndIf
        
    \EndWhile

    \State \hspace{-\fboxsep}\colorbox{lightgray!50!white}{\parbox{\dimexpr\linewidth-13pt-\algorithmicindent-2\fboxsep}{
    $K_0 \gets \arg\min_{\lfloor K/2 \rfloor < k < K-1} \|x^{k+1} - x^k\|$
    \hfill \textcolor{black}{\textit{Averaging} }
    }}
    
    \State $\hat{x} = \frac{1}{K_0 + 1} \sum_{k=0}^{K_0} z^k$
        \State \hspace{-\fboxsep}\colorbox{lightgray!50!white}{\parbox{\dimexpr\linewidth-13pt-\algorithmicindent-2\fboxsep}{
    $\nabla_\Theta L \mathrel{+}= JFB(\mathcal{L},T_\Theta^{\mathrm{\method}},\hat{x}, x^\star)$
    \hfill \textcolor{black}{\textit{Gradient update} }
    }}
\EndFor

\State \hspace{-\fboxsep}\colorbox{backbg}{\parbox{\dimexpr\linewidth-\algorithmicindent-1pt}{$\Theta \gets Adam(\Theta,\nabla_\Theta \mathcal{L})$ \hfill \textbf{Backward pass} }}

\EndFor

\end{algorithmic}
\rule{\linewidth}{0.8pt}
\end{minipage}

\vspace*{-0.7cm}
\end{wrapfigure}

where the model parameters $\Theta=\{\theta,\lambda\}$  are made learnable, yielding a fully data-driven but variationally consistent reconstruction model. We denote here $g_\theta$ instead of $g_\sigma$ to stress that the network weights are learned. The resulting equilibrium point $\bar x$ is defined as the fixed-point of the operator $T_\Theta^{\mathrm{\method}}$ and is taken as the reconstruction output.

A key distinction with standard DEQ formulations~\citep{Gilton2021DEQImaging,Zou2023ELDER} is the integration of inertia together with an adaptive restart strategy~\plaineqref{eq:restart} directly within the fixed-point computation. While inertia has already been used as an acceleration technique at inference time in PnP schemes~\citep{chow2024inertial,he2019plug,RISP, wu2024extrapolated}, we leverage it both during training and inference to improve the efficiency of the equilibrium solver.
During training, inertia improves the quality of the fixed-point estimate for a given number of iterations, providing a more accurate approximation of the equilibrium at a fixed computational budget and stabilizing training (see Figure \ref{fig:metric_vs_time} left). At inference time, it reduces the number of iterations required to reach a given fixed-point accuracy, thereby lowering the overall computational cost of the forward pass (see Figure \ref{fig:metric_vs_time} right and Figure \ref{fig:psnr_vs_iteration}). The adaptive restart mechanism further enhances stability by preventing the accumulation of excessive momentum, which can otherwise lead to oscillations or divergence, especially in nonconvex settings \citep{odono2015adaptive, RISP}.
Algorithm~\ref{algo:GM_RISP_DEQ} summarizes the inertial gradient-based scheme, while Algorithm~\ref{algo:RISP_DEQ_PROX} in Appendix~\ref{app:PGD} presents its proximal variant.

\begin{remark}
    Note that the averaging steps 12--13 of Algorithm~\ref{algo:GM_RISP_DEQ} are essential to ensure the convergence speed~\citep{li2023restarted}.
    They can be omitted without loss of performance by setting $\hat x=x^K$.
\end{remark}

\subsection{Convergence Analysis} \label{sec:convergence}

In this section, we establish the convergence of the algorithm to a fixed-point, in order to ensure the well-definedness of the method as defined in \citep{DEQ_Poisson}. We rely on the convergence results of~\citep{RISP} and make the following assumptions:

\begin{assumption}[Smoothness and Second-order smoothness]\label{ass:risp_1} The data-consistency is $L_f$-smooth, i.e. there exists $L_f > 0$ such that $\forall x,z \in \sR^d$, $ \|\nabla f(x,y) - \nabla f(z,y)\| \le L_f \|x - z\|$. Similarly, the regularization $g_\theta$ is $L_{g_\theta}$-smooth.

 The data-consistency has a $\rho_f$-Lipschitz Hessian, i.e., there exists $\rho_f > 0$ such that $\forall x,z \in \sR^d$, $\|\nabla^2 f(x,y) - \nabla^2 f(z,y)\| \le \rho_f \|x - z\|$. Similarly, the function $g_\theta$ has a $\rho_{g_\theta}$-Lipschitz Hessian.

\end{assumption}

\begin{remark}
    The smoothness assumption in~\ref{ass:risp_1} is standard in imaging inverse problems and is satisfied by many classical data-consistency terms, including linear inverse problems with additive Gaussian noise. The network requirement is mild and can be enforced using Lipschitz-continuous activations such as SoftPlus. Second-order smoothness also covers a broad class of linear inverse problems with additive white Gaussian noise, as well as certain nonlinear settings such as Rician noise models in MRI reconstruction \citep{gudbjartsson1995rician}. In practice, the network condition can be satisfied by combining gradient-step denoisers with Lipschitz activations; see~\citep[Appendix~G.1]{RISP} for details.
\end{remark} 

\begin{proposition}[{\citep[Theorems~1--2]{RISP}}] \label{prop:cv_DEQRISP}
Under Assumptions~\ref{ass:risp_1}, Algorithm~\ref{algo:GM_RISP_DEQ}, 
satisfy after at most $n$ iterations,
$\|\nabla F(\hat z)\| = \mathcal{O}(n^{-4/7})$.
\end{proposition}

This convergence rate improves upon the classical $\mathcal{O}(n^{-1/2})$ rate of standard RED methods and yields a provable acceleration for \method in reaching equilibrium during both training and inference.

\begin{remark}
    Proposition \ref{prop:cv_DEQRISP} guarantees the vanishing of the gradient norm ($\nabla F(x_k) \to 0$), and any point satisfying $\nabla F(x^*) = 0$ corresponds to a fixed point of the operator dynamics. 
    However, these theoretical results do not formally show point-wise convergence of the sequence $(x_k)_k$ to a unique fixed point. 
    In non-convex optimization, fixed points are generally not unique, and the algorithm is thus guaranteed to converge toward a critical point rather than a unique global fixed point.
\end{remark}

\section{Numerical Experiments} \label{sec:experiments}
We evaluate \method along three complementary directions. 
First, we conduct an extensive study on MRI reconstruction and analyze its training stability compared to classical DEQs. 
Second, we extend our method to image inpainting and motion deblurring as representative linear inverse problems. 
Finally, we consider reconstruction under Rician noise, a non-linear inverse problem where certain assumptions from Section~\ref{sec:convergence} no longer hold, to assess the robustness of our approach beyond its theoretical framework.
We evaluate \method against a wide spectrum of baselines: the E2E model RAM \citep{RAM}, PnP models such as RED/RED-P~\citep{HuraultGSD} and RISP/RISP-P~\citep{RISP}, the DEQ-based framework ELDER~\citep{Zou2023ELDER}, diffusion-based models including DiffPIR~\citep{Diffpir}, DPS~\citep{DPS}, and DDRM~\citep{DDRM}, as well as unrolled architectures like MoDL~\citep{aggarwal_modl_2019} and VarNet~\citep{hammernik_learning_2017} (see Appendices~\ref{app:DiffPIR}--\ref{app:unrolling_MRI} for details on baseline architectures and training). Note that specific baselines are omitted when they are not applicable to the forward model or fail to yield competitive results on the task considered.
Implementation and training details for \method are given in Appendices~\ref{app:implementation details} and~\ref{app:training_DEQ_gen}, while Appendix~\ref{app:hyperparameters} details hyperparameter choices for all baseline methods. We notably use the GS-Denoiser  architecture \plaineqref{eq:GS} for all the models but RAM to ensure a fair comparison. Codes are available at \url{https://github.com/antoninclerc/i-DEQ.git}

\subsection{Compressed sensing MRI} \label{sec:MRI}

We first evaluate \method on a compressed sensing MRI reconstruction task with an acceleration factor of $\times 8$, using the fastMRI dataset~\citep{knoll2020fastmri} restricted to the \texttt{knee-singlecoil} subset and extracting middle slices from 3D volumes. 
Gaussian noise with level $\sigma = 1/255$ is added in k-space, and the data are split into 100 training images, 10 validation images (used to tune hyperparameters and verify overfitting during training), and 20 test images (used to evaluate final performance).
In this setting, the data-consistency term satisfies Assumptions~\ref{ass:risp_1} under AWGN. The problem formulation  is detailed in~\ref{app:additional_MRI}.

\begin{figure}[t!]
\centering
\begin{tabular}{@{}c@{\hskip 1mm}c@{\hskip 1mm}c@{\hskip 1mm}c@{\hskip 1mm}c@{}}

% Ligne 1
\begin{tikzpicture}[spy using outlines={rectangle, orange, magnification=1.5, size=1.5cm}]
\node (img) {\includegraphics[width=0.18\linewidth, trim=0.7cm 0.5cm 0.7cm 1cm, clip]{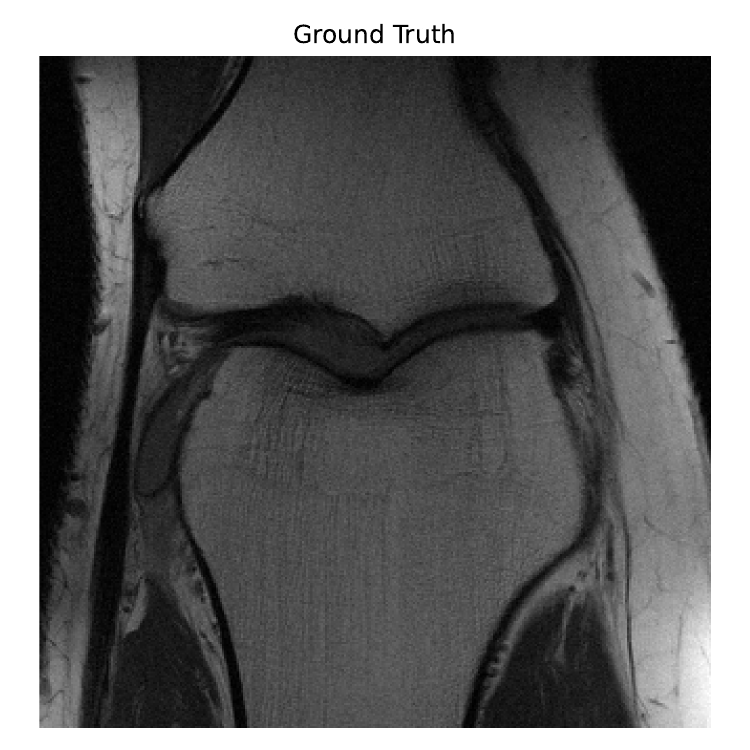}};
\spy on ($(img.center)+(-0.75,-0.25)$)
    in node [anchor=south east] at ($(img.south east)+(-0.05,0.05)$);
\end{tikzpicture}
&
\begin{tikzpicture}[spy using outlines={rectangle, orange, magnification=1.5, size=1.5cm}]
\node (img) {\includegraphics[width=0.18\linewidth, ]{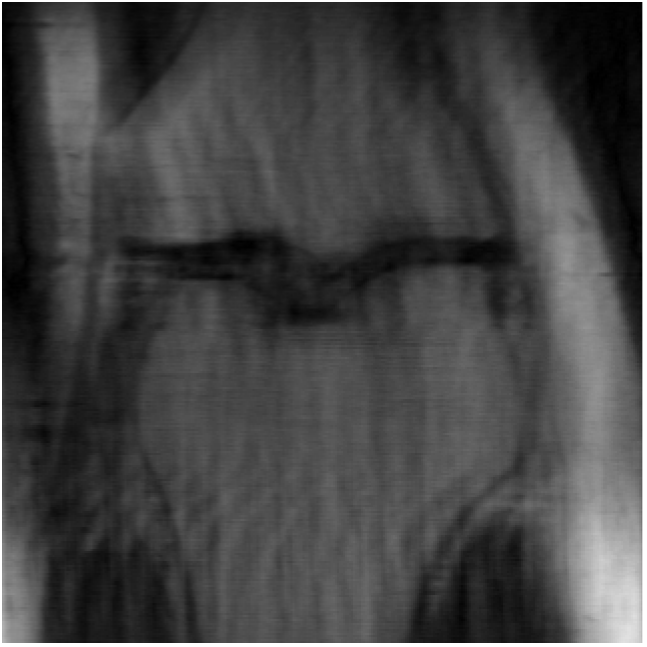}};
\spy on ($(img.center)+(-0.75,-0.25)$)
    in node [anchor=south east] at ($(img.south east)+(-0.05,0.05)$);
\end{tikzpicture}
&
\begin{tikzpicture}[spy using outlines={rectangle, orange, magnification=1.5, size=1.5cm}]
\node (img) {\includegraphics[width=0.18\linewidth, trim=1cm 0.5cm 1cm 1.5cm, clip]{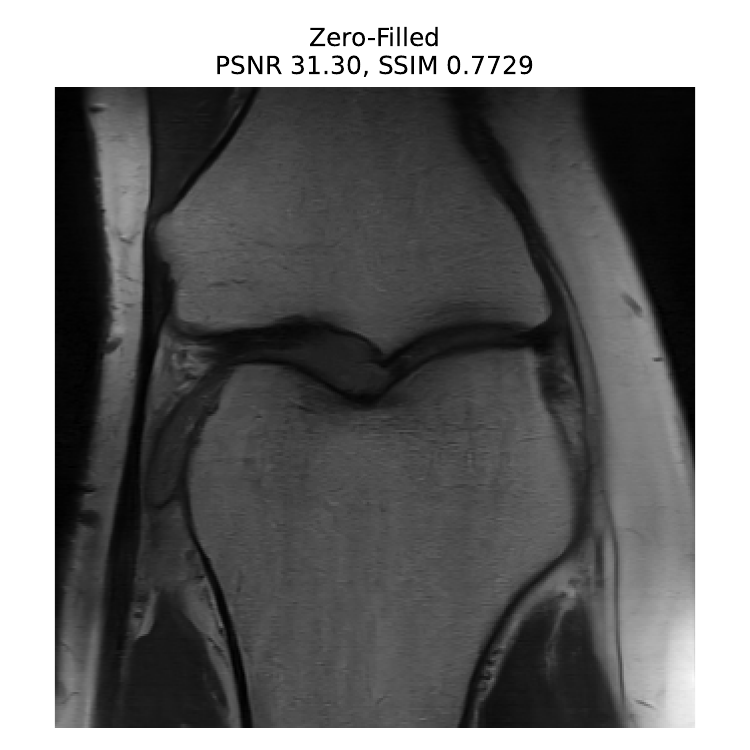}};
\spy on ($(img.center)+(-0.75,-0.25)$)
    in node [anchor=south east] at ($(img.south east)+(-0.05,0.05)$);
\end{tikzpicture}
&
\begin{tikzpicture}[spy using outlines={rectangle, orange, magnification=1.5, size=1.5cm}]
\node (img) {\includegraphics[width=0.18\linewidth]{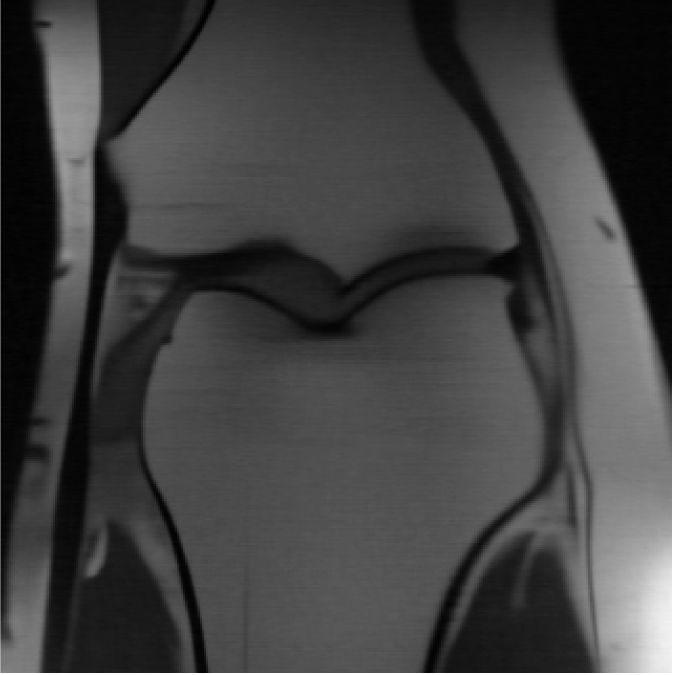}};
\spy on ($(img.center)+(-0.75,-0.25)$)
    in node [anchor=south east] at ($(img.south east)+(-0.05,0.05)$);
\end{tikzpicture}
&
\begin{tikzpicture}[spy using outlines={rectangle, orange, magnification=1.5, size=1.5cm}]
\node (img) {\includegraphics[width=0.18\linewidth, trim=0.8cm 0.5cm 0.8cm 1.1cm, clip]{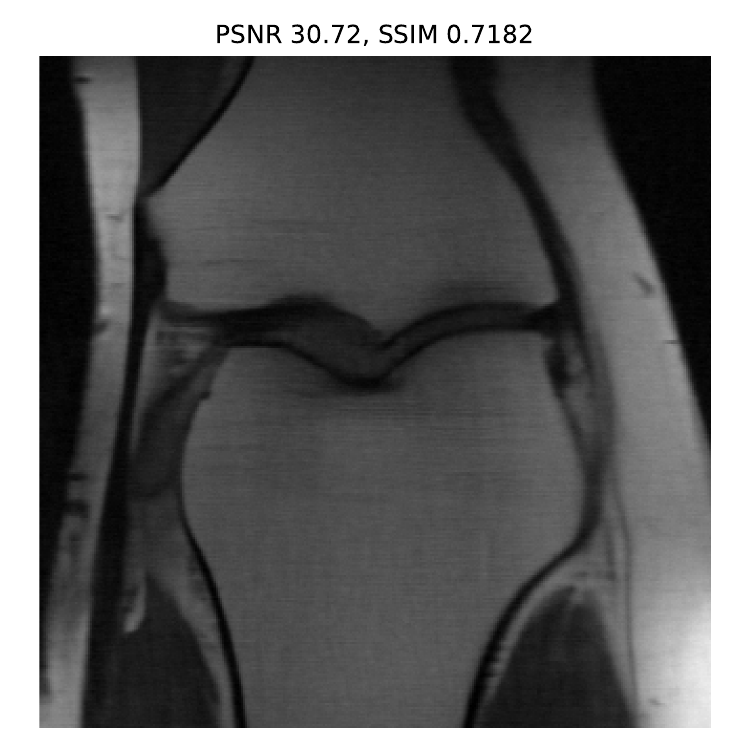}};
\spy on ($(img.center)+(-0.75,-0.25)$)
    in node [anchor=south east] at ($(img.south east)+(-0.05,0.05)$);
\end{tikzpicture}
\\
(a) Ground Truth & (b) Zero-filled & (c) RAM & (d) RISP & (e) DiffPIR \\
PSNR/SSIM & 23.86/0.59 & 31.30/0.77 & 31.26/0.74 & 30.72/0.72 \\

% Ligne 2`
\begin{tikzpicture}[spy using outlines={rectangle, orange, magnification=1.5, size=1.5cm}]
\node (img) {\includegraphics[width=0.18\linewidth]{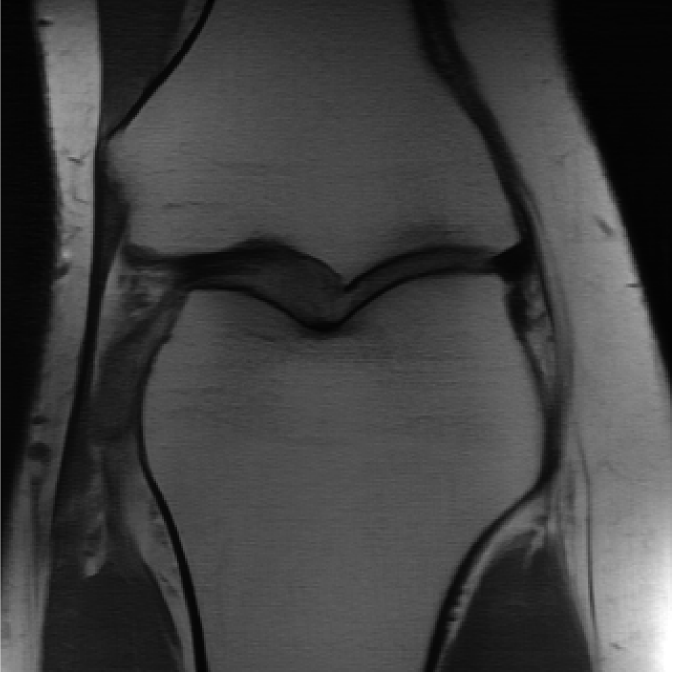}};
\spy on ($(img.center)+(-0.75,-0.25)$)
    in node [anchor=south east] at ($(img.south east)+(-0.05,0.05)$);
\end{tikzpicture}
&'
\begin{tikzpicture}[spy using outlines={rectangle, orange, magnification=1.5, size=1.5cm}]
\node (img) {\includegraphics[width=0.18\linewidth]{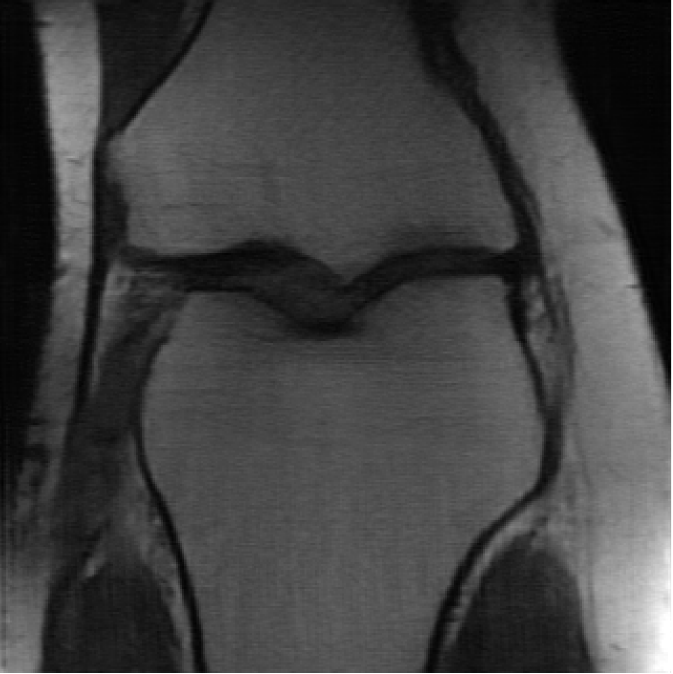}};
\spy on ($(img.center)+(-0.75,-0.25)$)
    in node [anchor=south east] at ($(img.south east)+(-0.05,0.05)$);
\end{tikzpicture}
&
\begin{tikzpicture}[spy using outlines={rectangle, orange, magnification=1.5, size=1.5cm}]
\node (img) {\includegraphics[width=0.18\linewidth]{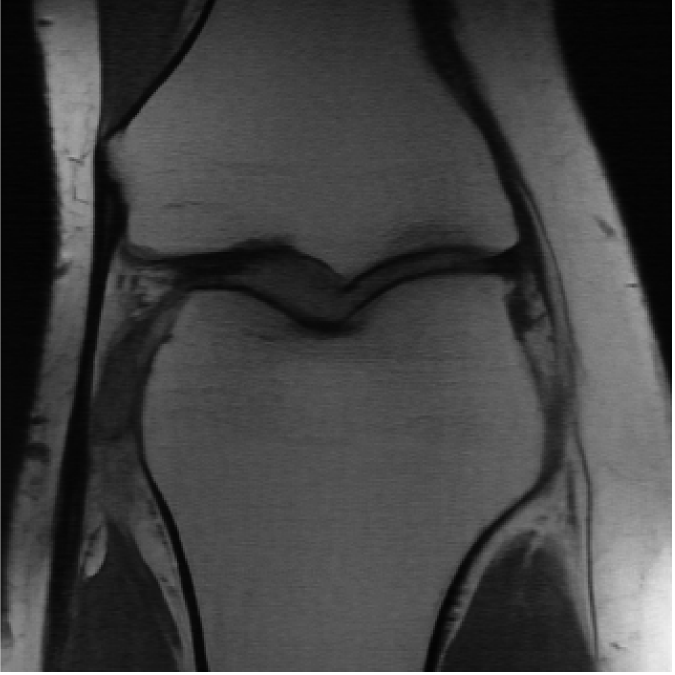}};
\spy on ($(img.center)+(-0.75,-0.25)$)
    in node [anchor=south east] at ($(img.south east)+(-0.05,0.05)$);
\end{tikzpicture}
&
\begin{tikzpicture}[spy using outlines={rectangle, orange, magnification=1.5, size=1.5cm}]
\node (img) {\includegraphics[width=0.18\linewidth]{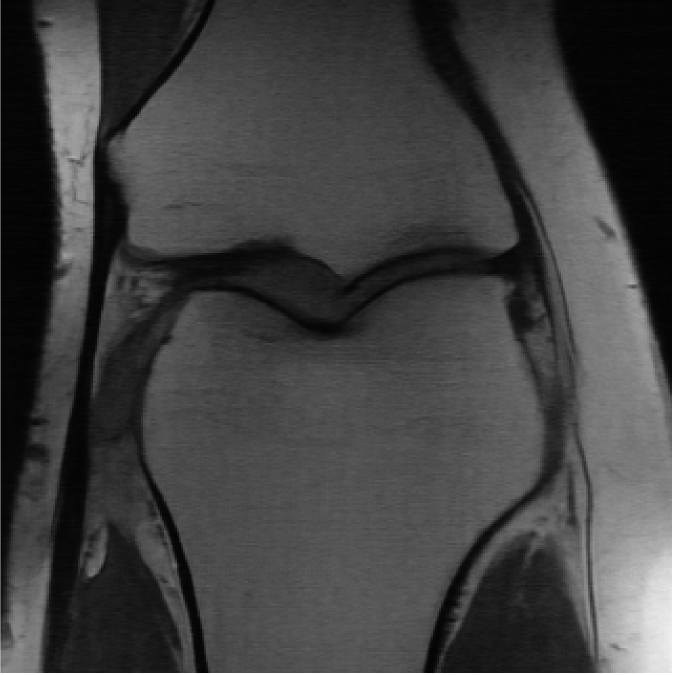}};
\spy on ($(img.center)+(-0.75,-0.25)$)
    in node [anchor=south east] at ($(img.south east)+(-0.05,0.05)$);
\end{tikzpicture}
&
\begin{tikzpicture}[spy using outlines={rectangle, orange, magnification=1.5, size=1.5cm}]
\node (img) {\includegraphics[width=0.18\linewidth, trim=1cm 0.5cm 1cm 1.5cm, clip]{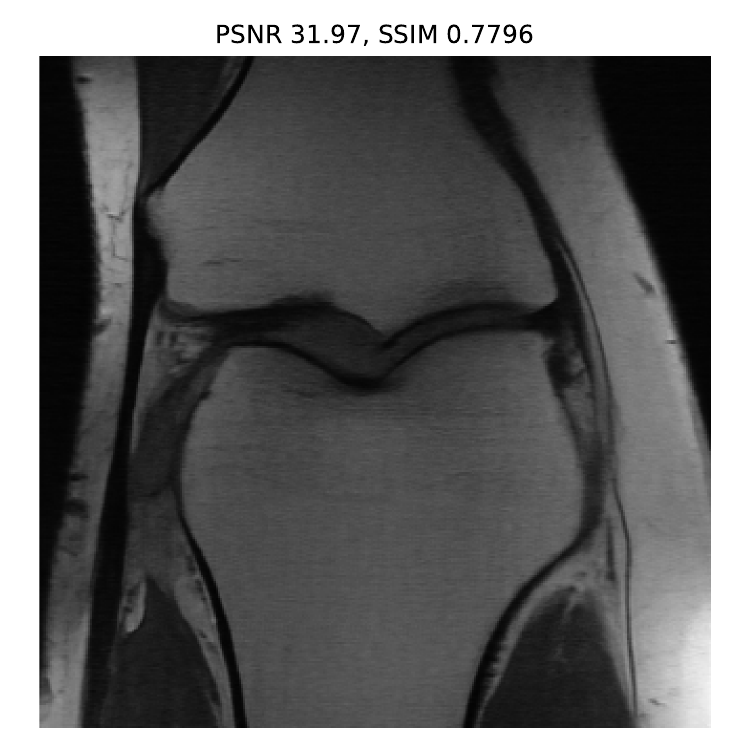}};
\spy on ($(img.center)+(-0.75,-0.25)$)
    in node [anchor=south east] at ($(img.south east)+(-0.05,0.05)$);
\end{tikzpicture}
\\
(f) MoDL & (g) VarNet & (h) ELDER & (i) \method & (j) \method-P \\
30.15/0.75 & 28.91/0.72 & 31.71/0.78  & 31.76/0.78 & 31.97/0.79
\end{tabular}
\caption{Compressed sensing MRI reconstruction with an acceleration factor of 8 using various iterative methods. DEQ approaches better preserve fine structural details than competing methods.
}
\label{fig:results_mri}
%\vspace{-0.5cm}
\end{figure}

\paragraph{MRI reconstruction}
Quantitative results in terms of PSNR and SSIM are reported in Table~\ref{table:psnr_mri}.
ELDER, \method and \method-P achieve the highest PSNR and SSIM among all evaluated methods. In contrast, unrolling approaches exhibit signs of overfitting during training, which leads to degraded generalization performance on the test set, but remain the fastest reconstruction techniques tested. RAM achieves better PSNR than unrolling approaches for a comparable computational cost, but falls short in term of PSNR compared to \method-P (-0.22dB).
DiffPIR achieves lower performance in term of PSNR (-0.57dB), comparable to RISP, which we attribute to the use of the light architecture of GS-DRUNet compared to the original work.
DDRM achieves performance comparable to DEQ-based approaches but inference takes more than twice as long as \method, and does not minimize any explicit energy.
Visual comparisons in \twoFigref{fig:results_mri}{fig:mri_x8_noise_1} show that DEQ reconstructions better preserve fine details and exhibit less oversmoothing than Plug-and-Play methods. 
\method is faster than ELDER at inference as shown in \Figref{fig:metric_vs_time} (right). When trained with a reduced number of fixed-point iterations, 60 seconds are sufficient to reach a PSNR of 28.19dB with \method compared to more than 300 seconds required by ELDER to achieve a similar reconstruction quality. This corresponds to a five times speed-up without loss of performance.

\paragraph{Anderson acceleration} 
Anderson Acceleration (AA) is a widely used approach for accelerating DEQ at inference \citep{Gilton2021DEQImaging}. We integrate it within a RED framework; further details are provided in Appendix~\ref{app:Anderson}.
We report the best-performing configuration, using a memory of five iterations without restart, for which no convergence guarantee exists in the non-convex setting. This configuration achieves a lower PSNR than \method (-0.33) and requires approximately 100s compared to 78s for \method to reach its peak performance, which is substantially longer due to the higher computational cost of the AA step. Overall, these factors make AA less suitable than \method.

\begin{table}[t!]
\centering\setlength{\tabcolsep}{4pt}\renewcommand{\arraystretch}{0.9}
\resizebox{\textwidth}{!}{
\begin{tabular}{l c c c c c c c c c c c c c c}
\toprule
& & E2E & \multicolumn{4}{c}{PnP} & \multicolumn{2}{c}{Diffusion} & \multicolumn{2}{c}{Unrolling} & \multicolumn{4}{c}{DEQ} \\
\cmidrule(lr){3-3} \cmidrule(lr){4-7} \cmidrule(lr){8-9} \cmidrule(lr){10-11} \cmidrule(lr){12-15}
Metric 
& Zero-filled & RAM & RED & RED-P & RISP & RISP-P
& DiffPIR & DDRM
& MoDL & VarNet 
& ELDER & Anderson & \method & \method-P  \\
\midrule
PSNR & 23.91 & 28.06 & 27.31 & 27.62 & 27.75 & 27.95 & 27.71 & 28.18 & 27.52 & 26.85 & 28.20 & 27.95 &\underline{28.21} & \textbf{28.28} \\
SSIM & 0.49 & \underline{0.62} & 0.54 & 0.56 & 0.56 & 0.56 & 0.56 & \underline{0.62} & \underline{0.62} & 0.60 & \textbf{0.63} & 0.49 & \textbf{0.63} & \textbf{0.63} \\
Time (s) & - & 1& 163 & 163 & 132 & 111 & 71 & 167 & \underline{0.6} & \textbf{0.5} & 319 & 100 & 70 & 78 \\
\bottomrule
\end{tabular}
}
\caption{Quantitative comparison in terms of PSNR and SSIM for MRI reconstruction on the test set across PnP, diffusion-based, unrolled optimization, and DEQ methods.  Best and second best results are respectively displayed in bold and underlined.  DEQ models achieve the best overall performance, with consistent improvements over both PnP and unrolled baselines.}
\label{table:psnr_mri}
%\vspace{-0.5cm}
\end{table}

\paragraph{Training stability} \label{sec:training_stability}
As illustrated in Figure \ref{fig:metric_vs_time} (left), the number of iterations used to approximate the fixed-point during the forward pass plays a critical role in DEQ training. With only $100$ iterations, a classical DEQ quickly becomes unstable and fails to converge. Increasing the budget to $200$ iterations improves its behavior, but does not fully prevent instability, as evidenced by a sharp PSNR drop occurring around $10^4$s (highlighted by a vertical marker in the figure).
In contrast, thanks to inertia, \method enables a fast, reliable estimation of fixed-points, resulting in stable training even with only $100$ fixed-point iterations. This increased stability comes at the cost of a slower PSNR improvement per epoch. However, this is offset by a significantly lower computational cost per iteration compared to classical DEQs, which require a backtracking strategy to avoid divergence during fixed-point computation. We also trained classical DEQ without backtracking to reduce the computational burden, but this resulted in degraded performance, with a loss of more than 1 dB. 

\paragraph{Robustness to parameter initialization} 
We evaluate the robustness of \method and ELDER with respect to the hyperparameter initialization ($\lambda$ and $\sigma$). Remarkably, \method trained with $200$ iterations exhibits high stability, remaining insensitive to initial hyperparameter values provided they enable early convergence. In contrast, classical DEQ is highly sensitive to initialization; deviating from grid-search-optimized values leads to a performance drop of over $1$\,dB in PSNR (Appendix~\ref{app:training_DEQ}). 
Regarding the momentum parameter $\alpha$, although smaller values theoretically speed up convergence, they impair training stability and yield longer inference times in practice. We thus fix $\alpha=0.2$, which achieves an effective compromise between convergence speed, robustness, and reconstruction fidelity (see Appendix~\ref{app:ablation_alpha}). Overall, these properties make \method more robust and virtually tuning-free, providing a major practical edge over standard DEQs.

\subsection{Inpainting} \label{sec:inpainting}

\begin{wraptable}{r}{0.54\textwidth}
\vspace{-0.35cm}
\centering\setlength{\tabcolsep}{4pt}\renewcommand{\arraystretch}{0.9}
\resizebox{0.55\textwidth}{!}{%
\begin{tabular}{l c c c c c c}
\toprule
& \multicolumn{3}{c}{Inpainting $\sigma=1/255$}
& \multicolumn{3}{c}{Inpainting $\sigma=5/255$} \\
\cmidrule(lr){2-4}
\cmidrule(lr){5-7}
Method
& PSNR & SSIM & Time
& PSNR & SSIM & Time \\
\midrule
Observation
& 8.47 & 0.13 & -
& 8.55 & 0.13 & - \\
RAM
& 33.37 & 0.91 & \bf{1}
& 32.56 & 0.89 & \bf{1} \\
RISP
& 33.87 & 0.93 & \underline{23}
& 33.29 & 0.90 & 42 \\
DiffPIR
& 34.63 & 0.93 & 37
& 33.45 & 0.90 & 34\\
DPS
& 30.76 & 0.84 & 95
& 28.96 & 0.78 & 95 \\
DDRM
& 34.84 & 0.93 & 45
& 33.86 & \underline{0.91} & 45 \\
ELDER
& 34.90 & 0.94 & 67
& \multicolumn{3}{c}{\textit{*training fails*}} \\
\method
& \underline{35.11} & \underline{0.95} & 25
& \underline{33.97} & \bf{0.92} & 44 \\
\method-P
& \bf{35.95} & \bf{0.96} & 30
& \bf{34.01} & \bf{0.92} & \underline{32}\\
\bottomrule
\end{tabular}
}
\caption{Quantitative comparison for image inpainting under multiple noise levels.} 
\label{table:inpainting}
\vspace{-0.45cm}
\end{wraptable}

We then consider an inpainting task where 50\% of the pixels are randomly removed, and Gaussian noise with levels $\sigma \in \{1, 5\}/255$ is added. Evaluations are conducted on a BSDS500 split \citep{amfm_pami2011}, consisting of 100 training images, 10 validation images, and 20 test images. A detailed formulation of the inverse problem and additionnal vidual results are provided in Appendix~\ref{app:additional_inpaint}.

As shown in Table~\ref{table:inpainting}, DEQ architectures outperform RISP and RAM at a low noise level ($\sigma = 1/255$). Although DPS falls short quantitatively, other generative approaches, namely DDRM and DiffPIR, achieve competitive results. Overall, \method and \method-P deliver the best reconstructions, outperforming classical DEQ, which suffers from training instabilities, and yielding superior visual results, as shown in Figure~\ref{fig:inpainting_sigma_1}.
At $\sigma = 5/255$, \method becomes more unstable during training and eventually diverges after $35$ epochs. Nevertheless, selecting the best validation checkpoint via early stopping yields a distinct gain over RISP. Meanwhile, \method-P remains stable throughout training. Both \method and \method-P maintain their edge over RISP ($+2.08$\,dB), DPS ($+5.19$\,dB), DiffPIR ($+1.32$\,dB), and RAM ($+2.58$\,dB), while performing on par with DDRM ($+0.27$\,dB for \method, $+1.11$\,dB for \method-P). Furthermore, \method(-P) retains the key theoretical advantage of minimizing an explicit energy formulation compared to diffusion and RAM approaches, while offering superior qualitative results. Conversely, the ELDER formulation diverges from the initial epochs, even with PnP parameter initialization, failing to produce usable models. Visual comparisons in Figure~\ref{fig:results_inpainting_rician_deblur} confirm that \method(-P) preserves fine details better and mitigates blurring artifacts compared to other approaches. 

\subsection{Motion Deblurring}

Next, we benchmark \method on a motion deblurring task corrupted by Gaussian noise ($\sigma = 5/255$), following the same dataset split and experimental protocol as in the inpainting setup. The complete mathematical model and extra visual examples are available in Appendix~\ref{app:additional_deblur}.

As detailed in Table~\ref{table:deblurring}, \method delivers the highest overall reconstruction quality, outperforming Prox-RISP by $+1.24$\,dB in PSNR. Both ELDER and DiffPIR yield comparable reconstruction metrics.
Regarding computational efficiency, \method offers a better trade-off than ELDER, cutting inference time from $184$s to $82$s, nearly a two-fold speedup. Moreover, \method is less prone to training instabilities than ELDER, whose performance deteriorates after just 6 epochs. While DiffPIR runs faster at inference, it requires an extensive hyperparameter search to stay competitive and lacks an explicit energy minimization framework compared to DEQ approaches. Visual results in Figure~\ref{fig:results_inpainting_rician_deblur} confirm that \method preserves fine details better and mitigates blurring artifacts compared to RISP-P.

\subsection{Rician denoising}
Finally, we extend our evaluation to Rician denoising under noise levels $\sigma \in \{0.05, 0.1\}$. Unlike previous tasks, this non-linear problem induces a non-convex data-consistency term $f$, offering a challenging testbed for algorithm stability. Experiments rely on the same BSDS500 split as in Section~\ref{sec:inpainting} to maintain protocol consistency. Further problem formulation and additional visual output are detailed in Appendix~\ref{app:additional_rician}.

\begin{figure}[t!]
\centering
\begin{tabular}{@{}c@{\hskip 1mm}c@{\hskip 1mm}c@{\hskip 1mm}c@{}}

% Ligne 1
\begin{tikzpicture}[spy using outlines={rectangle, orange, magnification=1.5, size=1.5cm}]
\node (img) {\includegraphics[width=0.23\linewidth]{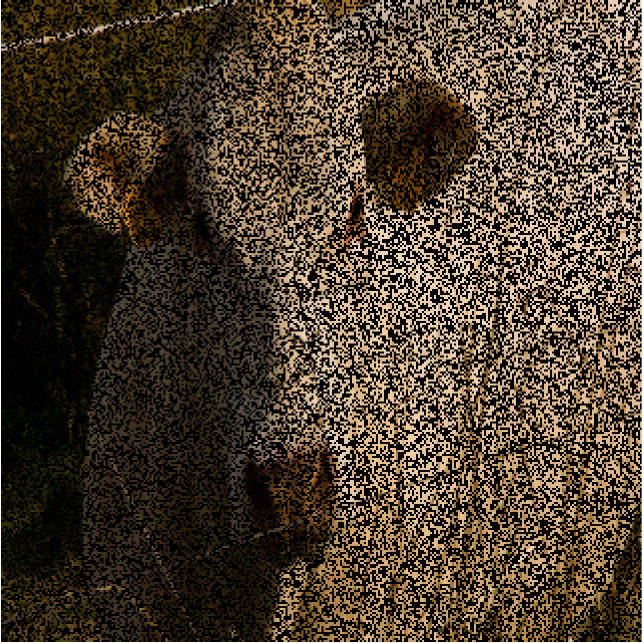}};
\spy on ($(img.center)+(0, 0)$)
    in node [anchor=south east] at ($(img.south east)+(-0.05,0.05)$);
\end{tikzpicture}
&
\begin{tikzpicture}[spy using outlines={rectangle, orange, magnification=1.5, size=1.5cm}]
\node (img) {\includegraphics[width=0.23\linewidth]{pics/test_reconstruction_inpaint_RISP.pdf}};
\spy on ($(img.center)+(0,0)$)
    in node [anchor=south east] at ($(img.south east)+(-0.05,0.05)$);
\end{tikzpicture}
&
\begin{tikzpicture}[spy using outlines={rectangle, orange, magnification=1.5, size=1.5cm}]
\node (img) {\includegraphics[width=0.23\linewidth]{pics/test_reconstruction_inpainting_ideq.pdf}};
\spy on ($(img.center)+(0, 0)$)
    in node [anchor=south east] at ($(img.south east)+(-0.05,0.05)$);
\end{tikzpicture}
&
\begin{tikzpicture}[spy using outlines={rectangle, orange, magnification=1.5, size=1.5cm}]
\node (img) {\includegraphics[width=0.23\linewidth]{pics/test_ground_truth_only_inpaint.pdf}};
\spy on ($(img.center)+(0, 0)$)
    in node [anchor=south east] at ($(img.south east)+(-0.05,0.05)$);
\end{tikzpicture}

\\
(a) Observation & (b) RISP & (c) \method & (d) Ground Truth \\
8.49/0.17   &  31.96/0.91 & 33.20/0.93 & PSNR/SSIM

\\

% Ligne 3
\begin{tikzpicture}[spy using outlines={rectangle, orange, magnification=1.5, size=1.5cm}]
\node (img) {\includegraphics[trim={1cm 0.5cm 1cm 1.4cm}, clip, width=0.23\linewidth]{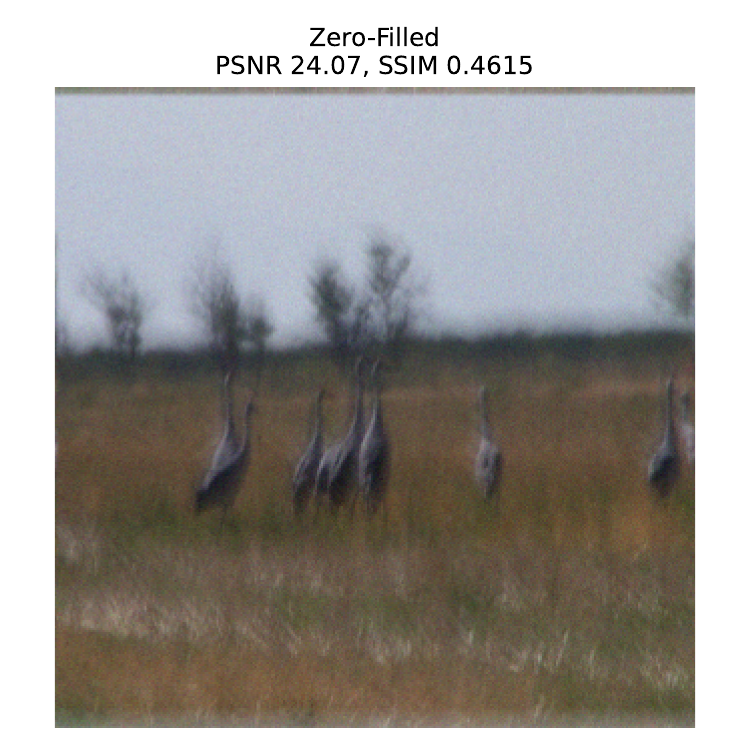}};
\spy on ($(img.center)+(-0.1,-0.53)$)
    in node [anchor=north west] at ($(img.north west)+(0.05,0)$);
\node[anchor=north east, draw=orange, line width=0.8pt, inner sep=0pt] at ($(img.north east)+(+0.05,0)$) {
    \includegraphics[trim={1cm 0.4cm 1cm 1.4cm}, clip, width=0.08\linewidth]{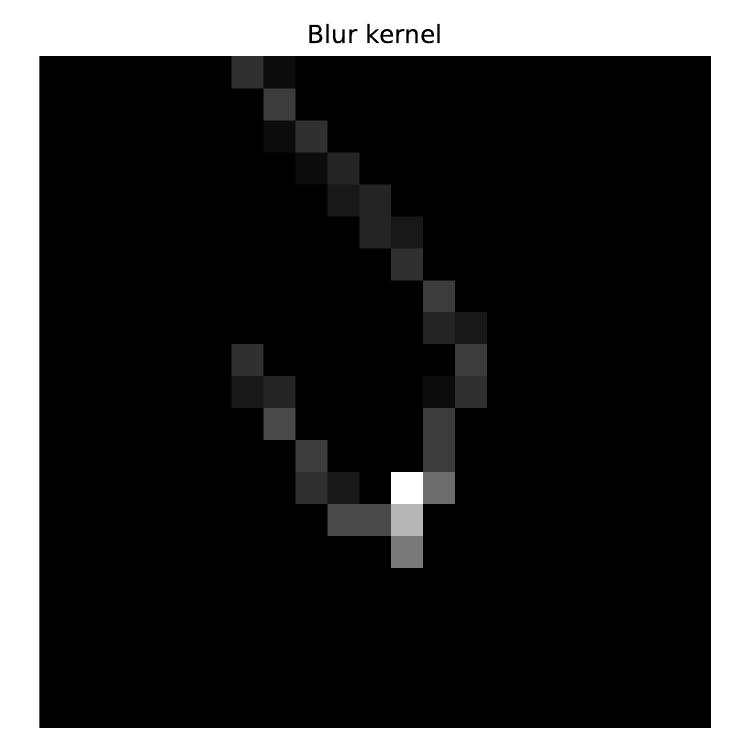}
};
\end{tikzpicture}
&
\begin{tikzpicture}[spy using outlines={rectangle, orange, magnification=1.5, size=1.5cm}]
\node (img) {\includegraphics[trim={1cm 0.5cm 1cm 1.4cm}, clip, width=0.23\linewidth]{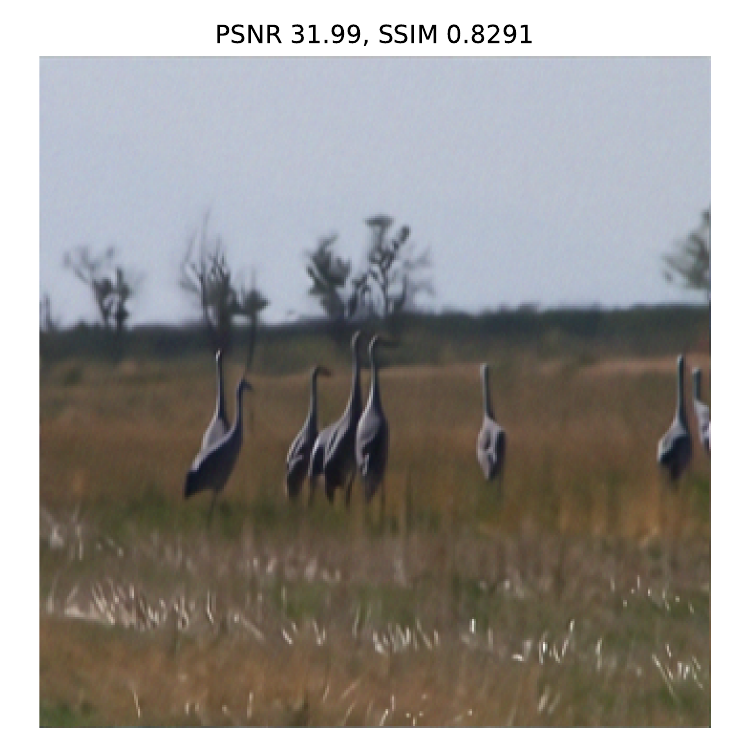}};
\spy on ($(img.center)+(-0.1,-0.45)$)
    in node [anchor=north west] at ($(img.north west)+(0.05,0)$);
\node[anchor=north east, draw=orange, line width=0.8pt, inner sep=0pt] at ($(img.north east)+(+0.05,0)$) {
    \includegraphics[trim={1cm 0.4cm 1cm 1.4cm}, clip, width=0.08\linewidth]{pics/test_blur_kernel.pdf}
};
\end{tikzpicture}
&
\begin{tikzpicture}[spy using outlines={rectangle, orange, magnification=1.5, size=1.5cm}]
\node (img) {\includegraphics[trim={1cm 0.5cm 1cm 1.4cm}, clip, width=0.23\linewidth]{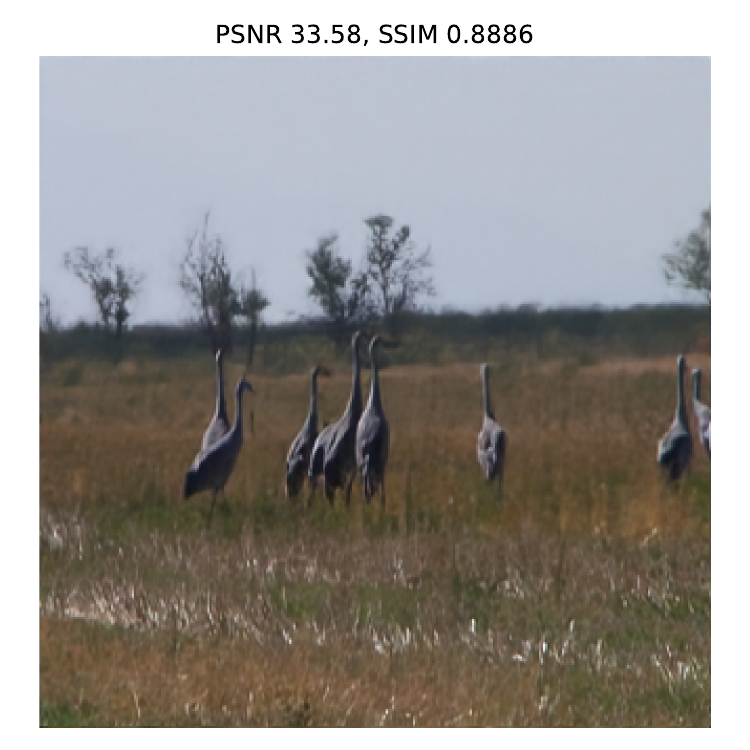}};
\spy on ($(img.center)+(-0.1,-0.45)$)
    in node [anchor=north west] at ($(img.north west)+(0.05,0)$);
\node[anchor=north east, draw=orange, line width=0.8pt, inner sep=0pt] at ($(img.north east)+(+0.05,0)$) {
    \includegraphics[trim={1cm 0.4cm 1cm 1.4cm}, clip, width=0.08\linewidth]{pics/test_blur_kernel.pdf}
};
\end{tikzpicture}
&
\begin{tikzpicture}[spy using outlines={rectangle, orange, magnification=1.5, size=1.5cm}]
\node (img) {\includegraphics[trim={1cm 0.5cm 1cm 1.4cm}, clip, width=0.23\linewidth]{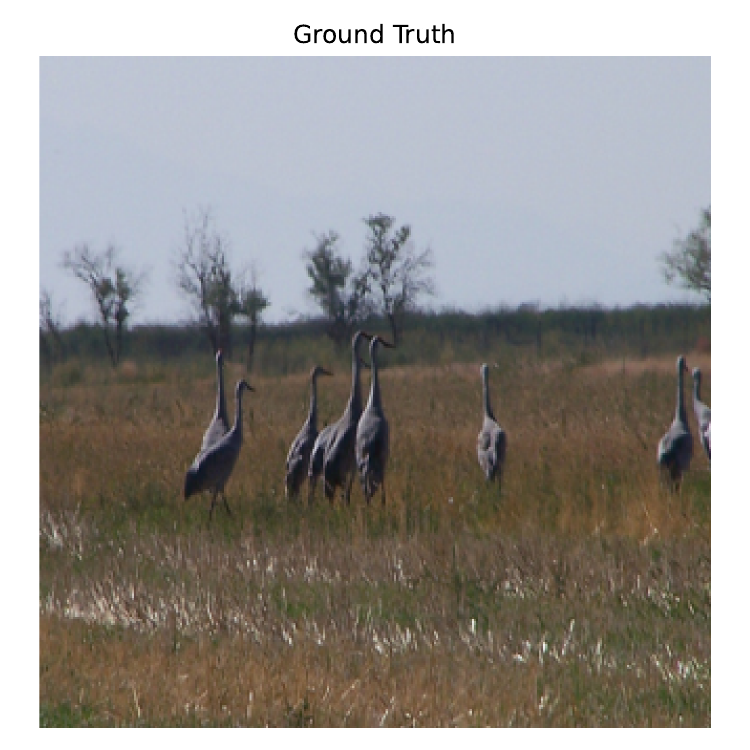}};

\spy on ($(img.center)+(-0.1,-0.45)$)
    in node [anchor=north west] at ($(img.north west)+(0.05,0)$);
\end{tikzpicture}
\\
(i) Observation & (j) RISP & (k) \method & (l) Ground Truth \\
24.07/0.46 & 31.99/0.83  & 33.58/0.89 & PSNR/SSIM \\

% Ligne 2
\begin{tikzpicture}[spy using outlines={rectangle, orange, magnification=1.5, size=1.5cm}]
\node (img) {\includegraphics[width=0.23\linewidth]{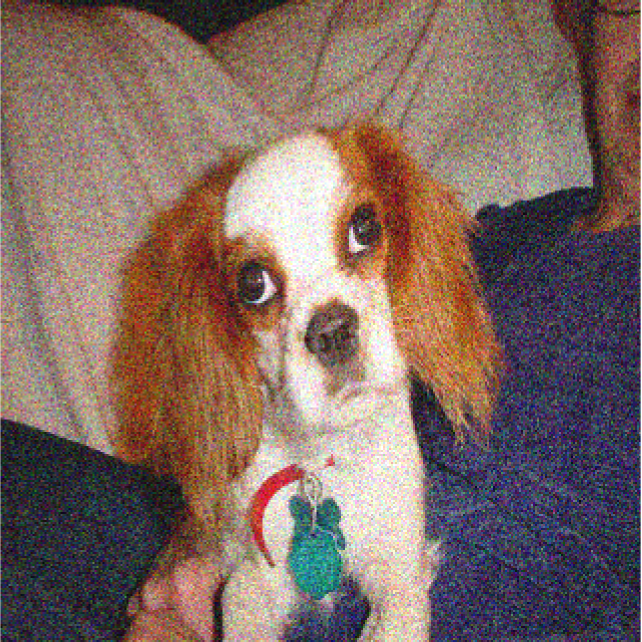}};
\spy on ($(img.center)+(0.85,-0.25)$)
    in node [anchor=north west] at ($(img.north west)+(0.05,0)$);
\end{tikzpicture}
&
\begin{tikzpicture}[spy using outlines={rectangle, orange, magnification=1.5, size=1.5cm}]
\node (img) {\includegraphics[width=0.23\linewidth]{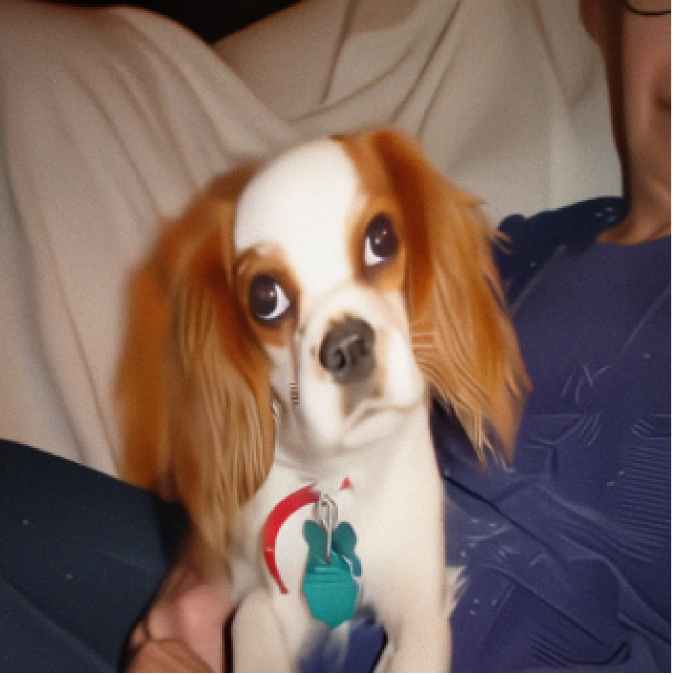}};
\spy on ($(img.center)+(0.9,-0.25)$)
    in node [anchor=north west] at ($(img.north west)+(0.05,0)$);
\end{tikzpicture}
&
\begin{tikzpicture}[spy using outlines={rectangle, orange, magnification=1.5, size=1.5cm}]
\node (img) {\includegraphics[width=0.23\linewidth]{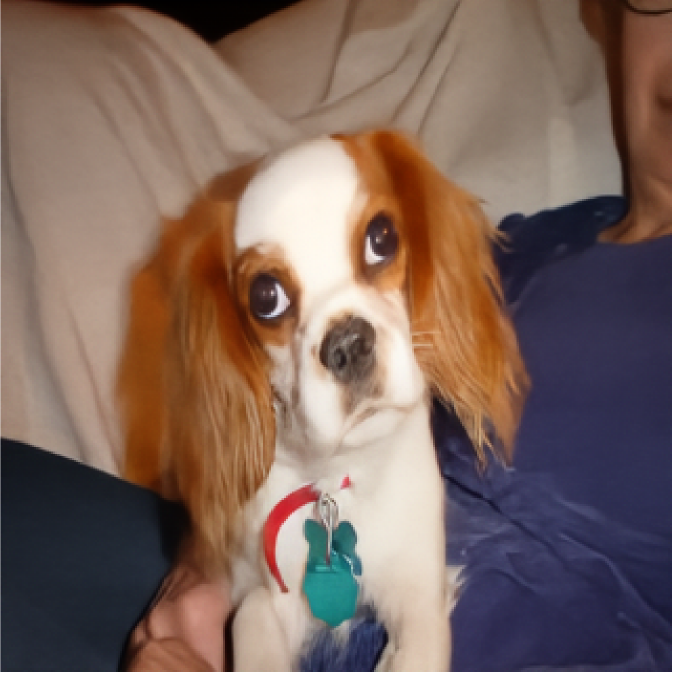}};
\spy on ($(img.center)+(0.9,-0.25)$)
    in node [anchor=north west] at ($(img.north west)+(0.05,0)$);
\end{tikzpicture}
&
\begin{tikzpicture}[spy using outlines={rectangle, orange, magnification=1.5, size=1.5cm}]
\node (img) {\includegraphics[width=0.23\linewidth]{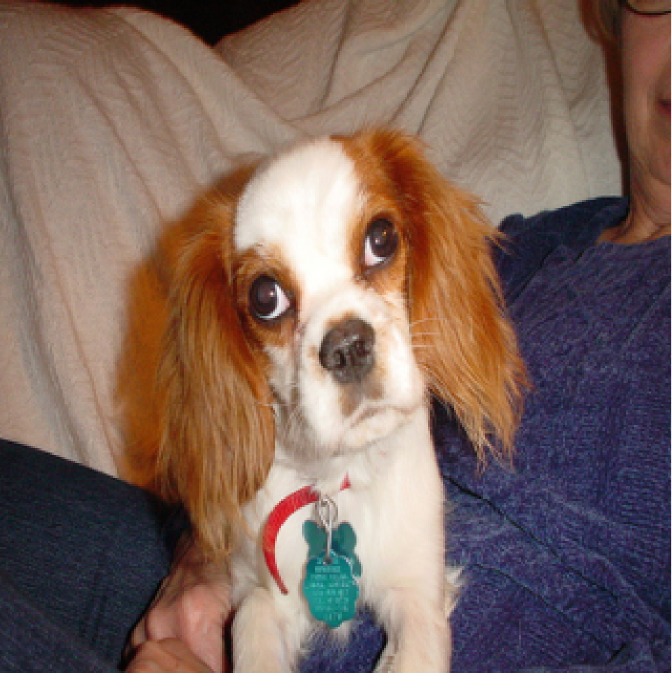}};
\spy on ($(img.center)+(0.9,-0.25)$)
    in node [anchor=north west] at ($(img.north west)+(0.05,0)$);

\end{tikzpicture}
\\

(e) Observation & (f) RISP & (g) \method & (h) Ground Truth \\
PSNR: 20.22/0.33 & PSNR: 30.38/0.77  & PSNR: 30.66/0.77 & PSNR/SSIM

\end{tabular}
\caption{The first row shows inpainting with 50\% of the pixels randomly removed and Gaussian noise with $\sigma = 5/255$. The second row shows motion deblurring with Gaussian noise of level $\sigma = 5/255$, and the third row shows Rician denoising with $\sigma = 0.1$. Across all three tasks, \method better preserves fine details while producing fewer artifacts than RISP.}

\label{fig:results_inpainting_rician_deblur}
\vspace{-0.3cm}
\end{figure}

As reported in Table~\ref{table:rician_denoising}, \method consistently outperforms all baselines across both noise regimes, achieving a PSNR gain over RISP of $+4.13$\,dB at $\sigma = 0.05$ and $+1.07$\,dB at $\sigma = 0.10$. Qualitative comparisons in Figure~\ref{fig:results_inpainting_rician_deblur} highlight how \method effectively suppresses noise and removes artifacts inherent to RISP, particularly in high-noise regimes. 
Conversely, ELDER systematically diverges during training due to the instability of its optimization trajectory; even when initialized with well-performing RISP parameters, it fails to converge. 
Note that generative-based baselines were omitted here, as they strictly assume additive Gaussian noise and are ill-suited for non-Gaussian statistics.

\begin{table}[t!]
\centering
% Premier tableau (Motion Deblurring)
\begin{minipage}[t!]{0.36\textwidth}
\centering
\setlength{\tabcolsep}{4pt}
\renewcommand{\arraystretch}{0.9}
\resizebox{\textwidth}{!}{%
\begin{tabular}{l c c c}
\toprule
& \multicolumn{3}{c}{Motion Deblurring} \\
\cmidrule(lr){2-4}
Method & PSNR & SSIM & Time \\
\midrule
Observation & 21.56 & 0.42 & - \\
RAM & 19.66 & 0.50 & \bf{1}\\
RISP-P       & 31.42 & 0.85 & 16.5 \\
DiffPIR     & \underline{32.54} & \underline{0.87} & {\underline{35}} \\
ELDER         & 32.47 & \underline{0.87} & 184 \\
\method-P    & {\bf 32.66} & {\bf 0.88} & 82 \\
\bottomrule
\end{tabular}%
}
\caption{Quantitative comparison for motion deblurring.}
\label{table:deblurring}
\end{minipage}% <--- Le % ici est crucial (pas de saut de ligne)
\hfill%
% Deuxième tableau (Rician Denoising)
\begin{minipage}[t!]{0.595\textwidth}
\centering
\setlength{\tabcolsep}{4pt}
\renewcommand{\arraystretch}{0.9}
\resizebox{\textwidth}{!}{%
\begin{tabular}{l c c c c c c}
\toprule
& \multicolumn{3}{c}{Rician $\sigma=0.05$}
& \multicolumn{3}{c}{Rician $\sigma=0.1$} \\
\cmidrule(lr){2-4}
\cmidrule(lr){5-7}
Method
& PSNR & SSIM & Time
& PSNR & SSIM & Time \\
\midrule
Observation
& 26.06 & 0.56 & -
& 20.19 & 0.33 & - \\
RISP
& \underline{31.33} & \underline{0.84} & \underline{14}
& \underline{30.85} & \underline{0.81} & 25 \\
ELDER
& \multicolumn{3}{c}{\textit{*training fails*}}
& \multicolumn{3}{c}{\textit{*training fails*}} \\
\method
& {\bf 35.46} & {\bf 0.92} & 25
& {\bf 31.92} & {\bf 0.86} & 35 \\
\bottomrule
\end{tabular}%
}
\caption{Quantitative comparison for Rician denoising under multiple noise levels.}
\label{table:rician_denoising}
\end{minipage}
\vspace{-0.5cm}
\end{table}

\section{Limitations}
Although i-DEQ provides significant acceleration compared to classical Plug-and-Play and DEQ-based approaches, several limitations remain.
First, although i-DEQ exhibits improved stability compared to classical DEQ architectures, training can still become challenging in extreme regimes, such as very high noise levels or severely ill-conditioned inverse problems. Further investigations into the stability properties of implicit models and the design of dedicated training strategies could improve robustness in these challenging settings.

Second, the theoretical analysis of i-DEQ relies on smoothness assumptions on the data-fidelity term, including Lipschitz continuity of the gradient and Hessian regularity. These assumptions enable the derivation of convergence and differentiation guarantees but exclude some important non-smooth formulations, such as $\ell_1$ data fidelity terms. Extending i-DEQ to such settings would require integrating advances in nonsmooth implicit differentiation and generalized derivative frameworks.

\section{Conclusion}
We introduced \method, an accelerated and tuning-free DEQ model that leverages inertial dynamics within both training and inference. First, we established the well-posedness of the formulation under mild assumptions. Then, through extensive experiments on linear and nonlinear inverse problems, we demonstrated that \method achieves state-of-the-art performance among iterative reconstruction methods. In particular, inertia significantly accelerates convergence to equilibrium, enabling \method to outperform classical DEQ models while reducing computation time by a factor of 2 to~5. Empirical results indicate that inertia improves training stability, allowing \method to be trained in more challenging settings while reducing the need for careful tuning of initialization parameters.

\section*{Acknowledgements}
AC, BD and NP gratefully acknowledges funding from the French Research Agency through the ANR-23-CE45-002 Holibrain and ANR-25-CE19-5558 OPLA.
MR gratefully acknowledges financial support from the French Direction Générale de l’Armement.

\newpage

\section*{Broader impacts}
This work on image reconstruction using deep equilibrium models may have several societal implications. First, reconstruction capabilities could be misused for malicious purposes such as the creation or enhancement of deepfakes, enabling deceptive practices including impersonation and manipulation of visual content. Second, \method, as a deep learning-based approach, when applied to sensitive domains such as medical imaging (e.g., MRI reconstruction), may inherit and amplify biases present in the training data. Moreover, as with other deep learning methods, training such models incurs non-negligible environmental costs due to energy consumption and carbon emissions, although reduced inference time compared to classical DEQs and iterative methods may partially offset this impact by lowering deployment-time energy usage.

\bibliographystyle{plain}
\bibliography{bib}

\appendix

\newpage

\section{More details on related works}

\paragraph{PnP}
PnP methods are a class of iterative algorithms for solving inverse problems of the form~\plaineqref{eq:inv_prob}, 
% combining a data-consistency term 
with a prior $g$ encoded by a denoiser. 
They were originally introduced as a modification of proximal optimization schemes, in which the proximal operator associated with the regularization term is replaced by an off-the-shelf denoiser ~\citep{venkatakrishnan2013plug}, without requiring a  formulation of the underlying regularizer.
A key limitation of this approach is that the denoiser does not, in general, correspond to the proximal operator or the gradient of any explicit regularization functional, unless additional assumptions are imposed~\citep{combettes2020lipschitz, hertrich2021convolutional, HuraultGSD, pesquet_learning_2021}. As a result, the underlying objective function is not explicitly defined, which complicates the analysis and weakens  convergence guarantees~\citep{ebner2024plug,ryu2019plug,wei2024learning}.

\paragraph{RED}
An alternative is the Regularization by Denoising (RED) framework~\citep{romano2017red}, which introduces an explicit regularization induced by a denoiser. In this setting, a denoiser $\Dsf_\sigma : \sR^d \rightarrow \sR^d$, with parameter $\sigma > 0$ controlling the denoising strength, is used to construct a gradient-like quantity through the residual $x - \Dsf_\sigma(x)$, interpreted as a surrogate for the gradient of an explicit regularizer, provided that the plugged denoiser has a symmetric Jacobian, which is an assumption that is hard to satisfy. In any case, this formulation provides a structured optimization viewpoint, enabling a tractable theoretical analysis of convergence from an operator-theory standpoint (i.e. convergence of iterates without targeting an explicit functional). Such analysis can be carried out either under Lipschitz-type assumptions on the denoiser~\citep{Reehorst2019,romano2017red}, or via explicit parameterizations of the underlying regularization functional~\citep{cohen2021it, fermanian2022learned, HuraultGSD}.

\paragraph{Generative models} Diffusion models (DMs)~\citep{ho2020denoising, sohl2015deep, Song2019, song2021score} and flow-matching (FM) methods~\citep{albergo2023stochastic, lipman2023flow, liu2023flow} have recently emerged as powerful generative priors, often interpreted as iterative denoising processes.
Their integration for solving inverse problems typically separates the data-consistency and prior modeling~\citep{chung2022diffusion, jalal2021robust}.
Therefore, pre-trained DMs~\citep{Diffpir, garber2024image, graikos2022diffusion, mardani2024variational} or FM models~\citep{martin2025pnpflow} can be plugged into iterative reconstruction schemes, leading to impressive reconstruction quality, especially in highly ill-posed settings.
However, theoretical guarantees for such restoration approaches remain limited~\citep{de2022convergence}.
In particular, it is not clear that these schemes optimize a target objective function such as the MAP estimator.
Moreover, to our knowledge, no existing approach has been rigorously shown to generate samples from the posterior distribution of the inverse problem~\plaineqref{eq:inv_prob}.

\paragraph{Deep unrolled networks}
Unrolled optimization methods construct deep architectures by truncating an iterative optimization algorithm to a fixed number of steps, where each iteration is mapped to a learnable network block. These approaches typically alternate between learned modules and explicit data-consistency steps and are trained end-to-end for a given forward model. They have achieved competitive performance across a wide range of imaging inverse problems~\citep{chun2020momentumnet, galinier2020hybrid, gharbi_unrolled_2024, mehranian2020model, monga_algorithm_2020, wu2019computationally}. Their application is particularly prominent in MRI, with representative works including \citep{aggarwal_modl_2019, hammernik_learning_2017, pramanik_adapting_2023, schlemper_deep_2017, sriram_end--end_2020, vo_efficient_2026}.

Formally, given an iterative scheme
\begin{equation}
x^{(k+1)} = T_{\theta^k}(x^{(k)}), \quad k = 0, \dots, K-1,
\end{equation}
unrolling consists in interpreting the composition of $K$ iterations as a feedforward network and learning the parameters $\theta$ by backpropagation through the associated computational graph. While this yields exact gradients for the truncated optimization procedure, it requires storing all intermediate iterates as well as all intermediate operations of the forward pass to enable reverse-mode differentiation.

Due to memory constraints, the number of iterations $K$ is typically kept small, often in the range of 5 to 10 in applications~\citep{pramanik_adapting_2023}.
Furthermore, unlike other model-based deep learning paradigms, the parameters $\theta$ are not necessarily shared across iterations and may explicitly depend on the iteration index $k$. This design choice moves the model further away from a classical iterative optimization algorithm, preventing it from accurately approximating a fixed point of the underlying dynamical system and bringing it closer to a standard end-to-end feedforward architecture.
Several works attempt to alleviate this limitation using memory-efficient strategies such as patch-based training \citep{vo_efficient_2026} or checkpointing techniques \citep{gruslys2016memory}, but these approaches generally still operate under a restricted number of unrolled steps. 
When weight sharing is enforced ($\theta^k = \theta$), other approaches attempt to train the network with a small number of unrolled iterations and reuse the learned weights for a larger number of iterations at inference time; however, such strategies often lead to performance degradation as the number of iterations increases~\citep{Gilton2021DEQImaging}.

The robustness of unrolling architectures is often assessed empirically through numerical experiments~\citep{genzel2023robustness,le2024unfolded}. More recently, theoretical studies have investigated stability and convergence of learning-based methods for inverse problems, notably via operator-theoretic approaches and structured proximal networks~\citep{hasannasab2020parseval,pesquet_learning_2021}. Complementary analyses rely on optimization-based interpretations, including unrolled and regularized networks, as well as variational and Lipschitz-based frameworks to establish stability guarantees~\citep{combettes2020variational,combettes2020lipschitz,genzel2023robustness,li2020nett}.

\paragraph{Comparison of the approaches.} 
Table~\ref{table:paradigm_comparison} provides a comprehensive comparison across paradigms, supported by empirical results from the compressed sensing MRI reconstruction task in Section~\ref{sec:MRI}.

\begin{table}[t!]
\centering\setlength{\tabcolsep}{4pt}\renewcommand{\arraystretch}{0.9}
\resizebox{\textwidth}{!}{
\begin{tabular}{p{2.2cm} p{2.2cm} p{2.2cm} p{2.4cm} p{2.4cm} p{2.6cm}}
\toprule
& \textbf{PnP} (\textit{RISP}) & \textbf{Unrolling} (\textit{MoDL}) & \textbf{Diffusion Models} (\textit{DDRM}) & \textbf{DEQ} (\textit{ELDER}) & \textbf{i-DEQ(-P) (Ours)} \\
\midrule
\textbf{Expressiveness (\textit{PSNR})} & 
  Low: generic prior \newline (\textit{27.95~dB}) & 
  Medium: task-specific priors with limited depth (\textit{27.52~dB}) & 
  High: deep generative priors \newline (\textit{28.18~dB}) & 
  High: implicit equilibrium deep prior, generalizing unrolled networks \newline (\textit{28.20~dB}) & 
  \textbf{High}: preserves DEQ expressiveness \newline (\textbf{\textit{28.28~dB}}) \\
\midrule
\textbf{MAP estimation} & 
  Converges toward MAP under standard PnP assumptions & 
  Does not approximate MAP or any explicit objective & 
  Samples from posterior, no/weak convergence results & 
  Convergence proofs under same assumptions as PnP & 
  \textbf{Convergence proofs} toward critical points \\
\midrule
\textbf{Inference time (\textit{time in s})} & 
  Slow, can be accelerated \newline (\textit{111s}) & 
  Extremely fast at inference \newline (\textit{$<$1s}) & 
  Slow (heavy sequential evaluations per step) \newline (\textit{167s}) & 
  Long inference \newline (\textit{319s}) & 
  \textbf{2$\times$ to 5$\times$ faster} than classical DEQs \newline (\textbf{\textit{78s}}) \\
\midrule
\textbf{Training} & 
  Separate database (pre-trained denoiser) & 
  End-to-end, memory scales with number of steps & 
  Separate score/diffusion pre-training & 
  Jacobian-free, constant memory w.r.t. steps & 
  \textbf{Jacobian-free} \\
\midrule
\textbf{Train. stability} & 
  High (standard denoiser training) & 
  Medium (end-to-end, prone to overfitting) & 
  High (stable score-matching objective) & 
  Low (often diverges or stagnates after few epochs) & 
  \textbf{Medium} (substantially more stable training) \\
\midrule
\textbf{Data Needs} & 
  Many available off-the-shelf & 
  Requires vast data to avoid overfitting & 
  Pre-trained on large datasets, needs matching target distribution & 
  Low, no substantial overfitting & 
  \textbf{Low, no substantial overfitting} \\
\bottomrule
\end{tabular}
}
\caption{Conceptual and empirical comparison across inverse problem paradigms. Quantitative results in italic are given for compressed sensing MRI reconstruction task \ref{sec:MRI}.}
\label{table:paradigm_comparison}
\end{table}

DEQ approaches carry the major advantage of being trainable on small datasets while maintaining high expressiveness. Beyond enhanced training stability, i-DEQ offers a clear variational interpretation, yielding convergence results compared to diffusion models and allowing sampling via Langevin dynamics. Furthermore, DEQ operates with lightweight neural network backbones. While diffusion models are highly competitive for general computer vision tasks with abundant data and computational resources, i-DEQ is better suited for clinical deployment scenarios where data are scarce and high-end GPU hardware is not readily available.
\section{Extension to the Proximal Gradient Descent Scheme} \label{app:PGD}

The proposed \method framework naturally extends to PGD schemes while preserving the convergence guarantees stated in Proposition~\ref{prop:cv_DEQRISP}, as detailed at the end of this section. This variant follows the formulation introduced in~\citep{Zou2023ELDER}, which relies on a proximal interpretation of RED-type methods.

Within the RED framework, the proximal-gradient-type iteration can be written as
\begin{equation}\label{eq:red_prox}
\textbf{RED--P:} \quad 
x^{k+1}
= T_{\theta, \lambda, \tau}^{\mathrm{RED\text{-}P}}(x^k)
\coloneqq
\operatorname{prox}_{\tau f}\!\Big(
x^k - \tau \lambda \nabla g_{\theta}(x^k)
\Big).
\end{equation}
where the relationship between the denoiser $\Dsf_{\sigma}$ and the regularization functional $g_{\theta}$ is given in~\plaineqref{eq:GS}.

The RISP framework~\citep{RISP} extends this scheme by incorporating inertia and adaptive restart~\plaineqref{eq:restart}  to stabilize convergence. Using the inertial operator defined in~\plaineqref{eq:Momentum},
and following the same proximal structure, 
we define the learned proximal operator \method-P by
\begin{equation}
T_{\Theta}^{\mathrm{\method\text{-}P}}(x^k, x^{k-1})
\coloneqq
T_{\theta, \lambda, \tau}^{\mathrm{RED\text{-}P}}\!\big(T^\mathrm{M}_\alpha(x^k, x^{k-1})\big),
\end{equation}
where $\Theta = \{\theta, \lambda\}$ are the learned parameters.
This formulation directly extends the PGD-based RISP scheme to the \method setting, while reducing the need for manual hyperparameter tuning by learning them directly from data. 
The corresponding fixed-point defines the output of the proximal \method model. Algorithm~\ref{algo:RISP_DEQ_PROX} summarizes the procedure, where inertia and restart components are highlighted for clarity.
Note that fixing the inertia parameter $\alpha = 1$, we recover the iteration of ELDER~\citep{Zou2023ELDER}.

\definecolor{forbg}{RGB}{245, 235, 200}
\definecolor{backbg}{RGB}{220, 240, 220}

\begin{algorithm}[t!]
\caption{\method-P}
\label{algo:RISP_DEQ_PROX}
\setlength{\fboxsep}{1pt}
\begin{algorithmic}[1]
\Require Dataset $(Y,X^\star)$, $n_{ep}, K \in \sN$, $\lambda, \tau > 0$, $\alpha \in (0,1)$, $B > 0$

\For{$n = 0$ to $n_{ep}-1$}

\State \hspace{-\fboxsep}\colorbox{forbg}{\parbox{\dimexpr\linewidth-\algorithmicindent-1pt}{$\nabla_\Theta L\gets 0$\hfill\textbf{Forward pass}}}

\For{each $(y,x^\star) \in (Y,X^\star)$} 
    \State $k \gets 0$, \quad $x^{-1} = x^0 \in \sR^d$
    
    \While{$k \leq K$}
    
        \State \hspace{-\fboxsep}\colorbox{lightgray!50!white}{\parbox{\dimexpr\linewidth-27pt-\algorithmicindent-2\fboxsep}{
        $z^k = x^k + (1 - \alpha)(x^k - x^{k-1})$
        \hfill \textcolor{black}{\textit{Inertial step}}
        }}
        
        \State $x^{k+1} \gets \prox_{\tau f(\cdot , y)} \left( z^k - \tau \lambda \nabla g_\theta(z^k) \right)$
        
        \State $k \gets k + 1$
        
        \If{$k \sum_{i=0}^{k-1} \|x^{i+1} - x^i\|^2 > B^2$}
            \State \hspace{-\fboxsep}\colorbox{lightgray!50!white}{\parbox{\dimexpr\linewidth-40pt-\algorithmicindent-2\fboxsep}{
            $x^{k-1} \gets x^k$ %$x^{-1}, x^0 \gets x^k$, \quad $k \gets 0$
            \hfill \textcolor{black}{\textit{Restart}}
            }}
        \EndIf
        
    \EndWhile

    \State \hspace{-\fboxsep}\colorbox{lightgray!50!white}{\parbox{\dimexpr\linewidth-13pt-\algorithmicindent-2\fboxsep}{
    $K_0 \gets \arg\min_{\lfloor K/2 \rfloor < k < K-1} \|x^{k+1} - x^k\|$
    \hfill \textcolor{black}{\textit{Averaging step}}
    }}
    
    \State $\hat{x} = \frac{1}{K_0 + 1} \sum_{k=0}^{K_0} z^k$
        \State \hspace{-\fboxsep}\colorbox{lightgray!50!white}{\parbox{\dimexpr\linewidth-13pt-\algorithmicindent-2\fboxsep}{
    $\nabla_\Theta \mathcal{L} \mathrel{+}= JFB(\mathcal{L},T_\Theta^{\mathrm{\method-P}},\hat{x}, x^\star)$
    \hfill \textcolor{black}{\textit{Gradient update}}}}
    
\EndFor

\State \hspace{-\fboxsep}\colorbox{backbg}{\parbox{\dimexpr\linewidth-\algorithmicindent-1pt}{$\Theta \gets Adam(\Theta,\nabla_\Theta \mathcal{L})$ \hfill \textbf{Backward pass}}}

\EndFor

\end{algorithmic}
\end{algorithm}

\begin{assumption}[Convexity structure]\label{ass:risp_2}
The function $f$ is convex, and there exists $\nu \ge 0$ such that the regularizer $g$ is $\nu$-weakly convex, i.e., the function
$
x \mapsto g(x) + \frac{\nu}{2}\|x\|^2
$
is convex.
\end{assumption}

\begin{remark}
Assumption~\ref{ass:risp_2} is mild in the sense that convexity of $f$ is standard in linear inverse problems. Moreover, if $g$ has an $L_g$-Lipschitz continuous gradient, then it is automatically $\nu$-weakly convex with $\nu \le L_g$ \citep{zhou2018fencheldualitystrongconvexity}.
The Gradient-Step denoiser with non SoftPlus activation used for the compressed sensing MRI reconstruction task (Section~\ref{sec:MRI}) does not strictly satisfy Assumption~\ref{ass:risp_2}. However, we observe in practice that it does not hinder convergence nor the acceleration effects predicted by the theory, in agreement with~\citep{RISP}.
\end{remark}

\begin{proposition}[{\citep[Theorem~2]{RISP}}] \label{prop:cv_DEQRISP_P}
Under Assumptions~\ref{ass:risp_1}-\ref{ass:risp_2}, Algorithm~\ref{algo:RISP_DEQ_PROX}, 
satisfies after at most $n$ iterations,
$\|\nabla F(\hat z)\| = \mathcal{O}(n^{-4/7})$.
\end{proposition}

Proposition~\ref{prop:cv_DEQRISP_P} provides the same convergence rate than Proposition~\ref{prop:cv_DEQRISP}. Note that the use of a proximal step on the data-fidelity may slightly improve the restoration in practice as observed in previous works~\citep{SNORE}.

\section{More details on competing methods}

\subsection{Anderson Acceleration} \label{app:Anderson}
Anderson acceleration (AA) \citep{Anderson1965,WalkerNi2011Anderson} is a widely used technique for accelerating fixed-point iterations by combining information from several previous iterates. Given a fixed-point scheme
$
x_{k+1}=T(x_k),
$
AA uses the residuals $r_i=T(x_i)-x_i$ from the previous $m$ iterations to compute a linear combination of past iterates that aims to reduce the fixed-point residual. The memory parameter $m$ therefore controls the amount of historical information used for acceleration. 
This type of extrapolation is widely used to accelerate fixed-point schemes and has been successfully applied to DEQ models and inverse problems \citep{Gilton2021DEQImaging}.
Despite its practical effectiveness, the convergence properties of AA are more challenging to establish for non-convex optimization problems. In particular, minimizing the fixed-point residual does not necessarily imply a decrease of the underlying objective function. Recent theoretical results show that restarted AA applied to gradient schemes can satisfy local descent properties under suitable conditions, providing a theoretical motivation for restart strategies in non-convex settings \citep{Ouyang_Liu_Milzarek_2024, Saad_2025}.

To evaluate AA in our MRI reconstruction setting, we consider memory depths $m\in{5,10}$, both with and without restarts. The results are reported in Table~\ref{tab}.
Among the tested AA configurations, the best-performing one in terms of reconstruction quality is obtained with a memory of $m=5$ without restart. Nevertheless, AA achieves lower reconstruction quality than \method at convergence. Although restarted AA reduces the number of iterations required for convergence, it also leads to a lower PSNR. Moreover, AA requires approximately 150 iteration to estimate the fixed point at each epoch to avoid divergence while only $100$ are sufficient for \method. 
Overall, our inertial strategy provides a more favorable trade-off between convergence speed, computational cost, and reconstruction quality.

\begin{table}[t!]
\centering
\footnotesize
\renewcommand{\arraystretch}{0.95}
\setlength{\tabcolsep}{5pt}
\caption{Comparison between Anderson acceleration and \method on the MRI reconstruction task.}
\label{tab}
\begin{tabular}{l c c c c}
\toprule
Solver / Config. & Restart & PSNR (dB) & SSIM & Time (s) \\
\midrule
Anderson ($m=5$) & No & 27.95 & 0.49  & 100 \\
Anderson ($m=5$) & Yes & 27.31 & 0.49  & 40 \\
Anderson ($m=10$) & No & 27.95 & 0.49 & 111 \\
Anderson ($m=10$) & Yes & 27.80 & 0.49& 66 \\
\method (100) & Yes ($B=100$) & 28.21 & 0.63 & 70 \\
\method (200) & Yes ($B=100$) & 28.21 & 0.63 & 144 \\
\bottomrule
\end{tabular}
\end{table}

\begin{remark}
    We do not consider standard Broyden's method \cite{Gilton2021DEQImaging, Broyden1965} in our experiments, as maintaining a full approximate Jacobian is particularly costly for the high-dimensional problems considered here. While limited-memory variants alleviate the memory requirements, they remain subject to similar local convergence considerations in non-convex settings as AA.
\end{remark}

\subsection{DiffPIR Algorithm}
\label{app:DiffPIR}

DiffPIR \citep{Diffpir} is an iterative Plug-and-Play (PnP) framework that integrates pre-trained diffusion priors with proximal operators to solve inverse problems of the form $y = \mathcal{A}(x_0) + \rvn$, where $\rvn \sim \mathcal{N}(0, \sigma_y^2 \mathbf{I})$ denotes additive Gaussian noise with standard deviation $\sigma_y$.

At each reverse step $t$ associated with noise schedule level $\sigma_t$, DiffPIR alternates between three main operations:
\begin{enumerate}[leftmargin=*]
    \item \textbf{Denoising (Prior):} Estimating the clean signal $\hat{x}_0(x_t)$ from the current noisy state $x_t$ using a pre-trained denoiser $D_{\sigma_t}$:
    \begin{equation*}
        \hat{x}_0(x_t) = D_{\sigma_t}(x_t).
    \end{equation*}
    \item \textbf{Data Consistency (Proximal Step):} Updating the estimated clean signal via a proximal mapping with respect to the measurement likelihood:
    \begin{equation*}
        \hat{x}_0(x_t) \leftarrow \mathrm{prox}_{\gamma_t f}\left(\hat{x}_0(x_t)\right), \quad \text{with } \gamma_t = \frac{\sigma_t^2}{2\lambda \sigma_y^2},
    \end{equation*}
    where $f(x) = \frac{1}{2}\|\mathcal{A}(x) - y\|_2^2$, and $\lambda > 0$ balances the data-fidelity term against the diffusion prior.
    \item \textbf{Reverse Diffusion Step:} Sampling the next state $x_{t-1}$ using the updated $\hat{x}_0(x_t)$ via a stochastic reverse step controlled by a variance parameter $\eta \in [0, 1]$.
\end{enumerate}

In our experiments, the denoiser $D_{\sigma_t}$ is instantiated using GS-DRUNet \citep{hurault2022proximal}. Detailed hyperparameter choices ($\lambda$, $\eta$, and step counts) across benchmark tasks are summarized in Table~\ref{table:hyperparameters_DiffPIR}.

\subsection{Diffusion Posterior Sampling (DPS)} 
\label{app:DPS}

Denoising Diffusion Probabilistic Models (DDPMs) \citep{ho2020denoising} define a forward process that gradually corrupts data $x_0$ into Gaussian noise via a prescheduled noise variance sequence $\bar{\alpha}_t \in (0, 1]$. Diffusion Posterior Sampling (DPS) \citep{chung2022diffusion} leverages such a pre-trained unconditional diffusion model as an unsupervised prior to solve general inverse problems of the form $y = \mathcal{A}(x_0) + \rvn$. It guides the reverse sampling trajectories using an approximate measurement likelihood gradient $\nabla_{x_t} \log p(y \mid x_t)$.

Since the exact conditional distribution $p(y \mid x_t)$ is intractable, DPS utilizes Tweedie's formula to estimate the posterior mean of the clean signal $\hat{x}_0(x_t) = \mathbb{E}[x_0 \mid x_t]$ from the current noisy state $x_t$ via the score network $s_\theta$:
\begin{equation*}
\hat{x}_0(x_t) = \frac{1}{\sqrt{\bar{\alpha}_t}} \left( x_t + (1 - \bar{\alpha}_t) s_\theta(x_t, t) \right).
\end{equation*}

At each reverse step, a standard DDPM transition yields an intermediate iterate $x'_{t-1}$, which is subsequently updated via gradient descent backpropagating through both the measurement operator $\mathcal{A}$ and the score network $s_\theta$:
\begin{equation*}
x_{t-1} \leftarrow x'_{t-1} - \zeta_t \nabla_{x_t} \|\mathcal{A}(\hat{x}_0(x_t)) - y\|_2^2,
\end{equation*}
where $\zeta_t$ denotes a normalized step size enforcing sampling stability.

In our experiments, we use $200$ reverse steps for MRI reconstruction and $300$ steps for natural image experiments to maintain an inference time comparable with competing approaches. Further hyperparameter specifications are detailed in Table~\ref{table:hyperparameters_DPS}.

\subsection{Denoising Diffusion Restoration Models (DDRM)}
\label{app:DDRM}

Denoising Diffusion Restoration Models (DDRM) \citep{kawar2022denoising} provide an unsupervised framework for linear inverse problems of the form $y = \mH x_0 + \rvn$, where $\rvn \sim \mathcal{N}(0, \sigma_y^2 \mathbf{I})$ denotes additive Gaussian noise with standard deviation $\sigma_y$.

Instead of relying on gradient-based likelihood guidance, DDRM performs inference directly in the spectral domain of the linear operator $\mH$. Utilizing the Singular Value Decomposition (SVD) $\mH = \mU \bm{\Sigma} \mV^\top$, the state $x_t$ and measurement $y$ are transformed into orthogonal subspaces via:
\begin{equation*}
\bar{x}_t = \mV^\top x_t, \qquad \bar{y} = \bm{\Sigma}^\dagger \mU^\top y,
\end{equation*}
where $\bm{\Sigma}^\dagger$ denotes the Moore--Penrose pseudo-inverse.

This decomposition decouples the inverse problem into independent 1D scalar sub-problems for each spectral component $i \in \{1, \dots, n\}$. At each step $t$ of the reverse sampling process, DDRM updates the spectral coordinates $\bar{x}_t^{(i)}$ by conditioning the transition distribution $p_\theta^{(t)}(\bar{x}_t^{(i)} \mid x_{t+1}, y)$ on both the unconditional score-based estimate $\hat{x}_0(x_{t+1})$ and the singular value $s_i$. The iterate in the spatial domain is subsequently reconstructed as $x_t = \mV \bar{x}_t$. By enforcing measurement consistency directly in the spectral domain, DDRM avoids costly backpropagation through the score network during sampling, offering a distinct computational advantage over guidance-based methods.

In our experiments, we use $200$ reverse steps for MRI reconstruction and $300$ steps for natural image tasks to maintain an inference time comparable to competing approaches. To ensure a fair comparison, the denoiser is instantiated using the exact same GS-DRUNet architecture \citep{hurault2022proximal} as in \method and DiffPIR.

\subsection{Unrolling models} \label{app:unrolling_MRI}

We train two classical unrolled models for MRI reconstruction:

\begin{itemize}
    \item \textbf{MoDL} \citep{aggarwal_modl_2019}:
    $
    x^{k+1} = \mathcal{N}_{\theta}\big(\operatorname{prox}_{\tau f}(x^k)\big),
    $
    
    \item \textbf{VarNet} \citep{hammernik_learning_2017}:
    $
    x^{k+1} = x^k - \tau \nabla f(x^k,y) - \mathcal{N}_{\theta}(x^k).
    $
\end{itemize}

In both cases, the algorithms are unrolled for 6 iterations. The step size $\tau$ is learned, and the network $\mathcal{N}_{\theta}$ is shared across all iterations. Although the performance of these architectures has been surpassed by more recent approaches \citep{liu2023diik,lu2025dffki,pramanik_adapting_2023,wang2025dsdun}, their structure, closely related to classical optimization schemes, makes them strong and widely used baselines for comparison with MBDL methods.

We use a DRUNet architecture for $\mathcal{N}_{\theta}$ to ensure consistency with the other numerical experiments. Unlike MBDL schemes, unrolled approaches do not impose specific structural constraints on the denoiser. This choice also avoids the additional memory overhead associated with backpropagation through the GS denoiser.
For MoDL, we use a pretrained DRUNet, initially trained on a denoising task, with weights from \citep{tachella2025deepinverse}, similarly to DEQ models. For VarNet, better results were empirically obtained when using a randomly initialized network.

For both models, the learning rate is set to $10^{-5}$, and we use the Adam optimizer with default hyperparameters. The models are trained to minimize the mean squared error (MSE) loss. Training is performed for up to 500 epochs with early stopping based on a patience criterion of 25 epochs.

\paragraph{Data Scaling and Overfitting Bottlenecks} 
To evaluate larger data regimes, we retrained MoDL using a $5\times$ larger fastMRI training set (500 images). Although additional data improves its performance in terms of PSNR/SSIM from 27.52 dB / 0.62 to 28.05 dB / 0.64, MoDL still overfits after approximately 25 epochs and remains below \method trained with only 100 images (28.21 dB / 0.63). 
These results highlight the improved data efficiency and generalization ability of DEQ-based models compared to unrolled architectures in data-scarce settings.

\section{More details on the experiments}

\subsection{Implementation details} \label{app:implementation details}

We reused the network architectures from the DeepInverse library~\citep{tachella2025deepinverse} and reimplemented the accelerated algorithms. For a fair comparison, all competing PnP and DEQ baselines use the same GS-Denoiser architecture as \method. In contrast, unrolling-based methods rely on a DRUNet architecture instead of the GS-Denoiser due to memory constraints. All experiments were conducted on a server equipped with NVIDIA RTX A6000 GPUs (48 GB VRAM).

For RGB  problems, the denoiser architecture is designed to satisfy Assumptions~\ref{ass:risp_1}. It relies on SoftPlus activation functions to ensure the required smoothness properties and is initialized from pretrained weights of~\citep{hurault2022proximal}. 
For grayscale problems, no pretrained denoiser satisfying the constraints is available. We instead use ELU activation functions, which have a Lipschitz-continuous gradient but a non-Lipschitz Hessian, and initialize the models using pretrained weights from~\citep{HuraultGSD}. 

We rely on JFB~\citep{JFB} for training \method, as is standard for DEQ models~\citep{DEQ_Poisson,Gilton2021DEQImaging,Zou2023ELDER}. Training details and hyperparameter choices are provided in Appendix~\ref{app:training_DEQ_gen}. Inference with \method is done by estimating a fixed-point with a forward pass.

\subsection{More details on DEQ training}\label{app:training_DEQ_gen}
In this section, we detail how the DEQ models (ELDER and \method) were trained in this paper. For simplicity, we call in this section "DEQ" the model "ELDER" proposed by~\citep{Zou2023ELDER}, the state-of-the-art DEQ model for image restoration.

\subsubsection{Reminder on the Backward Pass for DEQ}

Training a DEQ model requires differentiating through a fixed-point solution. Let $\bar{x}(\theta, y)$ denote the fixed point of the operator $T_\theta$, defined implicitly by
\[
\bar{x}(\theta, y) = T_\theta(\bar{x}(\theta, y), y).
\]
Given a dataset $\{(y_i, x_i^\ast)\}_{i=1}^n$, the training objective is
\[
\mathcal{L}(\theta) = \sum_{i=1}^n \ell\big(\bar{x}(\theta, y_i), x_i^\ast\big),
\]
where $\ell$ is typically the Mean Squared Error (MSE).
By the chain rule,
\begin{equation}
\nabla_\theta \ell(\bar{x}(\theta, y_i), x_i^\ast)
=
\left(\frac{\partial \bar{x}}{\partial \theta}(\theta, y_i)\right)^\top
\frac{\partial \ell}{\partial x}(\bar{x}(\theta, y_i), x_i^\ast).
\end{equation}

Backpropagating through the full forward solver is impractical due to the large number of iterations required for convergence and the resulting memory cost. Instead, differentiation in DEQ relies on the implicit differentiation of the fixed-point equation, which yields
\begin{equation}
\frac{\partial \bar{x}}{\partial \theta}(\theta, y)
=
\left(I - \frac{\partial T_\theta}{\partial x}(\bar{x}(\theta, y), y)\right)^{-\top}
\frac{\partial T_\theta}{\partial \theta}(\bar{x}(\theta, y); y).
\end{equation}

In this context, the inverse Jacobian is not computed explicitly in high dimensions, but is instead approximated using efficient schemes such as Jacobian-Free Backpropagation (JFB) \citep{JFB}. This corresponds to a first-order approximation of the inverse operator, leading to a practical gradient estimator \citep{bolte_one-step_2023}:
\begin{equation}
\nabla_\theta \mathcal{L}(\theta)
\approx
\sum_{i=1}^n
\left(\frac{\partial T_\theta}{\partial \theta}(\bar{x}(\theta, y_i); y_i)\right)^\top
\frac{\partial \ell(\bar{x}, x_i^\ast)}{\partial x}(\theta, y_i).
\end{equation}

However, backpropagation through DEQ models is known to be relatively unstable. A large body of literature therefore focuses on improving stability via implicit differentiation schemes, in particular by regularizing the Jacobian~\citep{bai2021stabilizing}, using alternative approximations to Jacobian-free backpropagation~\citep{blondel2022efficient, davy_restarted_2025,geng2022training}, or enforcing structural constraints that enable the computation of exact gradients~\citep{amos2021optnet,mccallum2025reversible}.

\subsubsection{Adaptation of Backward Pass for \method}

We now adapt the previous derivation for the operator $T^{\method}_\Theta(x_1,x_2)=T^{\method}(x_1,x_2, \Theta)$ to explicitly show the dependency in $\Theta=\{\theta,\lambda, \tau,\alpha\}$. We omit the observation $y$ in the operator and fixed point arguments for clarity. The fixed point $\bar{x}(\Theta)$ satisfies
\[
\bar{x}(\Theta) = T^{\method}(\bar{x}(\Theta), \bar{x}(\Theta), \Theta),
\]
Implicit differentiation gives
\begin{align*}
\frac{\partial \bar{x}(\Theta)}{\partial \Theta}
=&\;
\frac{\partial T^{\method}}{\partial x_1}
\frac{\partial \bar{x}}{\partial \Theta}(\Theta)
+
\frac{\partial T^{\method}}{\partial x_2}
\frac{\partial \bar{x}}{\partial \Theta}(\Theta)
+
\frac{\partial T^{\method}}{\partial \Theta},
\end{align*}
where all derivatives are evaluated at $(\bar{x}(\Theta),\bar{x}(\Theta),\Theta)$. Rearranging, we get
\begin{equation*}
\frac{\partial \bar{x}}{\partial \Theta}(\Theta)
=
\left(
I
-
\frac{\partial T^{\method}}{\partial x_1}
-
\frac{\partial T^{\method}}{\partial x_2}
\right)^{-1}
\frac{\partial T^{\method}}{\partial \Theta}.
\end{equation*}

Using the same JFB approximation leads to
\begin{equation}
\nabla_\Theta \mathcal{L}(\Theta)
\approx
\sum_{i=1}^n
\left(\frac{\partial T^{\method}}{\partial \Theta}(\bar{x}^{(i)}(\Theta), \bar{x}^{(i)}(\Theta), \Theta)\right)^\top
\frac{\partial \ell}{\partial x}(\bar{x}^{(i)}(\Theta), x_i^\ast).
\end{equation}

\begin{remark}
For the operator $T^{\method}$, the partial derivative with respect to the second argument vanishes at the exact fixed point. Likewise, derivatives with respect to $\tau$ and $\alpha$ are zero at convergence. In practice, the fixed point is only approximated, so these quantities are not exactly zero. However learning them has a negligible impact on the dynamics of \method\ during both training and inference, as shown in Section~\ref{app:training_DEQ}.
\end{remark}

\subsubsection{Training strategies} \label{app:training_DEQ}

All DEQ models are trained using a MSE loss $\mathcal{L}$ for up to 500 epochs, with early stopping (patience of 25 epochs) based on the PSNR measured on the validation set. The final model is selected as the checkpoint achieving the highest validation PSNR. Optimization is performed using the Adam optimizer with default parameters from the \textit{PyTorch} library~\citep{paszke2019pytorch} and a learning rate of $10^{-5}$.

Compared to classical DEQ approaches~\citep{Gilton2021DEQImaging,Zou2023ELDER}, we define fixed-point convergence in \method as the completion of a prescribed number of successive iterations without having reach the restart condition, in accordance with Algorithms~\ref{algo:GM_RISP_DEQ} and~\ref{algo:RISP_DEQ_PROX}, rather than relying on a convergence criterion based on the norm of successive iterates, i.e., $\|x^{k+1} - x^k\|$. This choice is motivated by the presence of inertia, which can induce larger oscillations between successive iterates even when the sequence is close to a fixed-point. As a result, the norm $\|x^{k+1} - x^k\|$ may not accurately reflect the quality of the fixed-point approximation, making such stopping criteria less reliable in this setting. For a fair comparison, classical DEQ approaches are trained using the same budget of iterations. Note that in practice, we still introduced the stopping criterion $\|x^{k+1} - x^k\| / \|x^k\|<\varepsilon$ with $\varepsilon=10^{-4}$ to avoid unnecessary computations once near convergence if the budget of iterations is too high; this criterion is occasionally triggered during the early epochs but is no longer used thereafter.

Table~\ref{table:hyperparameters_DEQ} reports the hyperparameter settings used to initialize DEQ models. The regularization parameter $\lambda$ and inertia parameter $\alpha$ are learned during training. The step size $\tau$ is either learned for \method or selected via a backtracking line search for ELDER. To distinguish the parameter $B$ values used during training from those applied at inference, we denote them as $B_{\text{train}}$ and $B_{\text{infe}}$, respectively. For most of the models, $B_{\text{train}} = B_{\text{infe}}$ ; in that case we write $B_{\text{train}} = B_{\text{infe}} = B$.
% The inertia parameter $\alpha$ is also learned during training.
The noise level $\sigma$ plays an important role at initialization, as the denoiser is pretrained over a range of noise levels. However, during DEQ training, the model is optimized for a single fixed noise level, making $\sigma$ effectively constant and having limited influence on the final model behavior. Finally,  the restart parameter $B$ is set to a lower value than in RISP (see Table~\ref{table:hyperparameters_PnP}), as this choice empirically improves training stability.

\begin{remark}
Learning the parameters $\alpha$, $\tau$, and $\lambda$ has empirically negligible impact on performance; they can be fixed without significant loss in accuracy.
\end{remark}

We trained several variants to analyze stability and the effect of the number of iterations:

\begin{itemize}
    \item \method (100 / 200): we consider that the fixed-point is reached when 100 (rep. 200) iterations have been realized without restart. In practice, a single restart is observed making the effective number of iterations used for the JFB step typically around 125 or 225. The initial hyperparameters are the optimal one obtained from a grid search on RISP.

    \item \method (100!): the fixed-point is estimated after exactly 100 iterations, regardless of any intermediate restart event.
    The initial hyperparameters are the optimal one obtained from a grid search on RISP.

    \item \method / \method-P / DEQ (arb.init.): variants where the initial parameters are not tuned to their optimal configuration; instead, they are deliberately set to smaller values called "arbitrary initialization" to study sensitivity to initialization.

    \item DEQ (100 / 150 / 200): the fixed-point is computed after 100, 150 or 200 iterations. No restart is used since this method does not include inertia or momentum mechanisms.
    The initial hyperparameters are the optimal one obtained from a grid search on RED.
\end{itemize}

The initial hyperparameter configurations are detailed in Table~\ref{table:hyperparameters_DEQ}.

\begin{table}[t!]
\centering \footnotesize
\begin{tabular}{lcccccccc}
\toprule
Inverse problem & Method & $\lambda$ & $\sigma$ & $\tau$ & $\alpha$ & $B$ & PSNR & SSIM \\
\midrule

\multirow{8}{*}{MRI $\times 8$}
& \method (100!)        & 0.65 & 0.03 & 0.5 & 0.2 & 100 & 28.02 & \textbf{0.63} \\
& \method (100)       & 0.65 & 0.03 & 0.5 & 0.2 & 100 & 28.16 & \textbf{0.63} \\
& \method (100)       & 0.65 & 0.03 & 0.5 & 0.2 & \shortstack{$B_{\text{train}} = 100$ \\ $B_{\text{infe}} = 50$} & 28.20 & \textbf{0.63} \\
& \method (200)       & 0.65 & 0.03 & 0.5 & 0.2 & 100 & 28.21 & \textbf{0.63} \\
& \method (200, arb.init.)   & 0.10 & 0.00 & 1.0 & 0.2 & 100 & 28.21 & \textbf{0.63} \\
& \method (100) & 0.8 & 0.02 & 1.0 & 0.2 & 100 & \textbf{28.28} & \textbf{0.63} \\
& \method-P (200, arb.init.) & 0.10 & 0.00 & 1.0 & 0.2 & 100 & 28.22 & \textbf{0.63} \\
& DEQ (100)            & 0.83 & 0.03 & 2.0   & /   & /   & 28.02 & 0.62 \\
& DEQ (200)            & 0.83 & 0.03 & 2.0   & /   & /   & 28.20 & \textbf{0.63} \\
& DEQ (200, arb.init.)       & 0.83 & 0.00 & 2.0   & /   & /   & 27.54   & 0.61     \\

\midrule

\multirow{2}{*}{Inpainting, $\sigma = 1/255$}
& \method (100) & 0.83 & 0.03 & 0.1 & 0.2 & 500 & \textbf{35.11} & \textbf{0.95} \\
& DEQ (150)      & 0.83 & 0.03 & 2.0 & /   & /   & 34.90 & 0.94 \\

\multirow{2}{*}{Inpainting, $\sigma = 5/255$}
& \method (100) & 0.83 & 0.03 & 0.1 & 0.2 & 500 & \textbf{33.97} & \textbf{0.92} \\
& \method (200) & 0.83 & 0.03 & 0.1 & 0.2 & 500 & 33.80 & \textbf{0.92} \\

\midrule
Rician, $\sigma = 0.05$ & \method (100) & 10.0 & 0.02 & 0.03 & 0.2 & 100 & \textbf{35.46} & \textbf{0.92} \\

\multirow{2}{*}{Rician, $\sigma = 0.1$}
& \method (100) & 6.0    & 0.02 & 0.03 & 0.2 & 300 & \textbf{31.92} & \textbf{0.86} \\
& \method (200) & 6.0    & 0.05 & 0.03 & 0.2 & 300 & 31.90 & \textbf{0.86} \\

\multirow{2}{*}{Deblurring}
& \method-P (100) ($\sigma = 5/255$) & 0.84    & 0.03 & 0.5 & 0.2 & 50 & \textbf{32.66} & \textbf{0.88} \\
& DEQ (150) & 0.84    & 0.03 & 0.5 & / & / & 32.48 & 0.87 \\
\bottomrule
\end{tabular}
\caption{Hyperparameter settings at the initialization of the DEQ models. If not specified, the parameter $B_{\text{train}}$ used during training is the same as $B_{\text{infe}}$ used during inference.}
\label{table:hyperparameters_DEQ}
\end{table}

\begin{figure}[ht]
    \centering
    \includegraphics[width=0.5\linewidth]{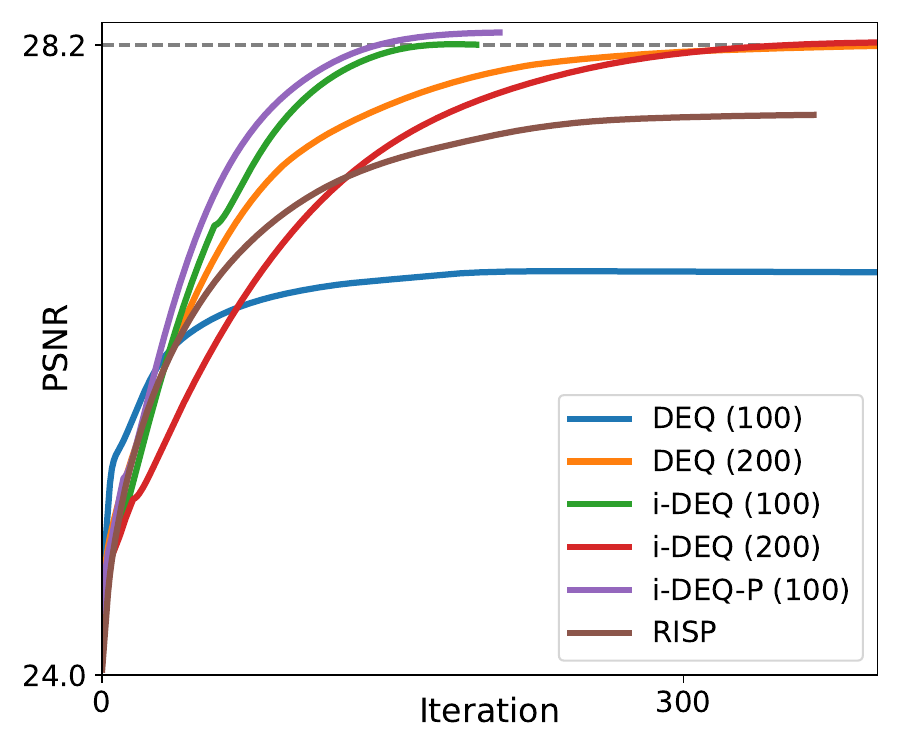}
    \caption{PSNR evolution versus inference iterations. Evaluation on 20 test images from the FastMRI dataset~\citep{knoll2020fastmri}. \method (trained with 100 and 200 fixed-point steps) is compared against classical DEQ baselines with matching iteration budgets, as well as RISP. \method (100) converges in $2$--$3\times$ fewer iterations than both 200-step variants (\method and DEQ). Although DEQ (200) exhibits a comparable convergence trajectory to \method (200), its per-iteration cost is substantially higher owing to backtracking operations.}
    \label{fig:psnr_vs_iteration}
\end{figure}

In \Figref{fig:psnr_vs_iteration}, we evaluate the evolution of PSNR with respect to the iteration count during reconstruction. Notably, even when trained with a restricted iteration budget, \method continues to improve reconstruction quality beyond its training depth. We thus evaluate both classical DEQ and \method at inference until reaching fixed-point convergence, defined by the stopping criterion $\|x^{(k+1)} - x^{(k)}\|_2 / \|x^{(k)}\|_2 < 10^{-4}$. Remarkably, \method (100) and \method-P (100) require approximately half as many iterations to converge as their 200-iteration counterparts, \method (200) and \method-P (200), without any drop in PSNR. While DEQ (200) requires a similar number of iterations to converge as \method (200), each iteration incurs a significantly higher computational cost due to backtracking line searches. Consequently, its total inference time is 2--5$\times$ higher than that of \method (see the last column of Table~\ref{table:hyperparameters_DEQ}). Thus, \method (100) and \method-P (100) offer the best trade-off between reconstruction quality and computational efficiency.

% \begin{remark}
% The temporary performance drop around 60 seconds does not signify a failure of convergence.
% This behavior is mathematically sound and perfectly consistent with convergence. Throughout this phase, the residual norm $\|z^{(k+1)} - z^{(k)}\|$ and the objective energy continue to decrease monotonically. This confirms that the solver remains numerically stable and does not diverge. \method minimizes a fixed-point energy formulation rather than directly maximizing the PSNR at every step. Due to the nature of the updates, the iterates can temporarily bypass the optimal visual reconstruction state (peak PSNR) while still strictly descending within the energy landscape, as observed by \citep{RISP, odono2015adaptive}. Our convergence Proposition (Proposition B.1) fully aligns with this behavior: it guarantees asymptotic convergence to a critical point with respect to this energy, without requiring a monotonic increase in PSNR. Therefore, the observed drop in PSNR while the energy decreases is not a mathematical aberration, but a direct consequence of the optimization trajectory.
%\end{remark}

\subsubsection{Sensitivity to the parameter B} 
As in previous section, we distinguish $B_{train}$ and $B_{infe}$ as potentially different restart criteria for training and inference.
We evaluate the sensitivity of i-DEQ across $B \in \{10, 50, 100, 200, 500\}$, considering that training convergence is reached when 150 forward iterations have been realized without restart (i.e. $150 \sum_{i=0}^{149} \|x^{i+1} - x^i\|^2 \leq B_{train}^2$).
%During training under a budget of 150 forward iterations without restart, 
For $B_{\text{train}} \in \{10, 50, 200\}$, restart events cause the effective number of inference iterations to fluctuate unpredictably across epochs on both training and validation sets. This variability impairs training stability and induces high variance in validation PSNR, complicating the monitoring of model stagnation or overfitting. In contrast, setting $B_{\text{train}} = 100$ provides a remarkably stable regime during validation and testing by consistently triggering a single restart around iteration 25, yielding a deterministic total of roughly 175 fixed-point iterations. Alternatively, taking $B = 500$ prevents any restart within the 150-iteration budget, effectively reducing the process to a standard fixed-iteration scheme. 
Crucially, when evaluation is governed by a fixed total iteration budget rather than an uninterrupted iteration threshold—i.e., when \method evaluates exactly 150 iterations per epoch—training dynamics remain essentially independent of $B_{\text{train}}$. Ultimately, training stability is governed by maintaining a predictable, steady iteration depth across epochs.

During inference, we run \method until a fixed point is reached, which we assess using the criterion $|x^{k+1} - x^k|^2 \leq 10^{-4}$. When trained with few iterations, the model may exhibit performance degradation at inference if $B_{\text{infe}}$ is chosen too large. Mathematically, while the objective energy and residual norm decrease monotonically, thus guaranteeing numerical stability and asymptotic convergence to a critical point (Proposition~\ref{prop:cv_DEQRISP} and \ref{prop:cv_DEQRISP_P}), this energy descent does not enforce a monotonic increase in PSNR, occasionally causing temporary drops \citep{RISP, odono2015adaptive}. In practice, taking $B_{\text{infe}}=50$ prevents this phenomenon. 
Specifically, for short trajectories, \method (100) achieves optimal performance when trained with $B_{\text{train}} = 100$ and evaluated with $B_{\text{infe}} = 50$ to prevent degradation, whereas \method-P reaches its optimum with $B_{\text{train}} = B_{\text{infe}} = 100$.

Overall, $B_{train}=B_{infe}=B = 100$ provides a robust default setting for both training and inference across variants, offering a tuning-free baseline without fine-grained adjustments.

\subsubsection{Influence of the momentum parameter $\alpha$}
\label{app:ablation_alpha}

We evaluate the effect of the momentum parameter $\alpha$ on the $\times 8$ MRI reconstruction task for $\alpha \in {0, 0.01, 0.1, 0.2, 0.5, 0.9}$. Note that when taking $\alpha=1$ we recover the classical DEQ formulation. We compare the different settings in terms of training stability, inference speed, and final reconstruction performance measured by PSNR.

\begin{figure}[ht]
    \centering
    \includegraphics[width=0.6\textwidth]{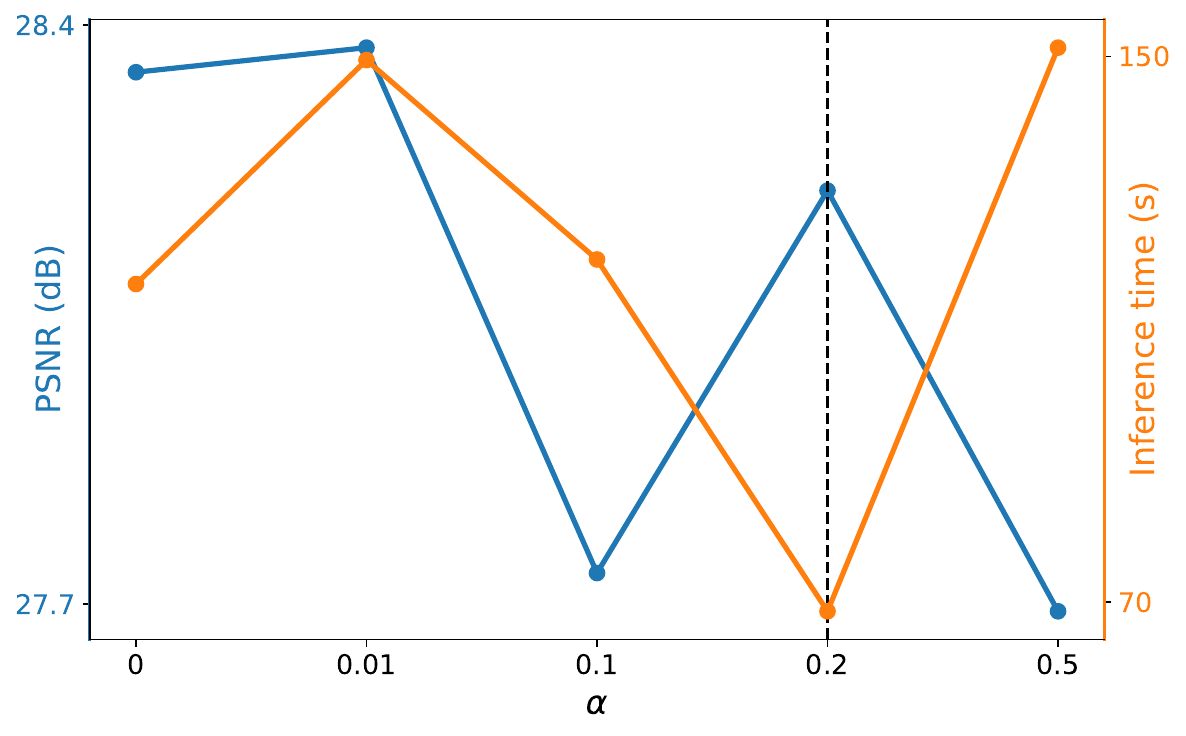}
    \caption{Ablation study on the effect of the momentum parameter $\alpha$. $\alpha=0.2$ provides the best compromise between inference speed and reconstruction performance.}
    \label{fig:ablation_alpha}
\end{figure}

For large values of $\alpha$, training becomes less stable or prone to overfitting. In particular, $\alpha=0.9$ leads to divergence, while $\alpha=0.5$ remains stable but overfits after a few epochs, resulting in a final PSNR of 27.70 dB. A similar degradation is observed for $\alpha=0.1$, with a final PSNR of 27.74 dB.
Conversely, smaller values of $\alpha$ lead to stable training but require substantially more iterations to reach convergence, with significant variations in the number of iterations across epochs. For $\alpha=0.01$, the final PSNR reaches 28.37 dB, while $\alpha=0$ achieves 28.34 dB. Although these values slightly outperform our default setting $\alpha=0.2$, the improvement of approximately 0.15 dB comes at a substantially higher inference cost. Moreover, $\alpha=0$ lies outside the assumptions of our convergence guarantees, which require $\alpha>0$.
\Figref{fig:ablation_alpha} summarizes the trade-off between inference time and final PSNR. Overall, $\alpha=0.2$ provides the best compromise between training stability, inference speed, and reconstruction performance, motivating its use as our default setting.

\subsection{Parameter setting }\label{app:hyperparameters}
\subsubsection{PnP approaches parameter settings}

Table \ref{table:hyperparameters_PnP} reports the hyperparameter settings of PnP methods selected via grid search for the different inverse problems. Inertia parameters are chosen according to \citep{RISP}.

\begin{table}[t!]
\centering

\begin{tabular}{lcccccc}
\toprule
Inverse problem & Method & $\lambda$ & $\sigma$ & $\tau$ & $\alpha$ & $B$ \\
\midrule

\multirow{4}{*}{MRI $\times 8$}
& RISP       & 0.65 & 0.03 & 0.5 & 0.2 & 5000 \\
& RISP--Prox & 0.80 & 0.02 & 1.0 & 0.2 & 5000 \\
& RED        & 0.80 & 0.05 & 0.5 & /   & /    \\
& RED--Prox  & 0.80 & 0.03 & 1.0 & /   & /    \\

\midrule

\multirow{2}{*}{Inpainting}
& RISP ($\sigma = 5/255$) & 0.83 & 0.03 & 0.1 & 0.2 & 5000 \\
& RISP ($\sigma = 1/255$) & 0.83 & 0.03 & 0.1 & 0.2 & 5000 \\

\midrule

\multirow{2}{*}{Rician}
& RISP ($\sigma = 0.05$) & 10   & 0.02 & 0.03 & 0.01 & 100 \\
& RISP ($\sigma = 0.1$)  & 3.6  & 0.03 & 0.03 & 0.01 & 100 \\
% & RISP-P ($\sigma = 0.05$) & 2320 & 0.03 & $4.3\times10^{-6}$ & 0.01 & 100 \\
% & RISP-P ($\sigma = 0.1$)  & 834  & 0.03 & $1.2\times10^{-5}$ & 0.01 & 100 \\

\midrule

Deblurring $\sigma = 5/255$
& RISP-P & 0.5 & 0.03 & 2.0 & 0.2 & 50 \\

\bottomrule
\end{tabular}

\caption{Hyperparameters used for the different inverse problems.}
\label{table:hyperparameters_PnP}

\end{table}

\subsubsection{DiffPIR model parameter settings}

Table \ref{table:hyperparameters_DiffPIR} reports the hyperparameter settings of DiffPIR selected via grid search for the different inverse problems. We use the GS-Denoiser \citep{HuraultGSD} architecture for all task for fair comparison with \method.

\begin{table}[t!]
\centering
\begin{tabular}{lcccc}
\toprule
Task & $\lambda$ & $\zeta$ & $N$ \\
\midrule

MRI & 1.0 & 0.19 & 200 \\

Inpainting ($\sigma = 1/255$) & 0.1 & 1.0 & 200 \\
Inpainting ($\sigma = 5/255$) & 0.63 & 1.0 & 200 \\

Deblurring ($\sigma = 5/255$) & 0.34 & 1.0 & 200 \\

\bottomrule
\end{tabular}
\caption{Hyperparameter settings for DiffPIR across different inverse problems.}
\label{table:hyperparameters_DiffPIR}
\end{table}

\subsubsection{DPS Hyperparameter Settings}
\label{app:DPS_hyperparameters}

Table~\ref{table:hyperparameters_DPS} summarizes the optimal hyperparameter configurations for DPS selected via grid search. In our implementation using the \texttt{deepinv.sampling.DPS} module \citep{tachella2025deepinverse}, $N$ denotes the total number of diffusion sampling steps (\texttt{num\_steps}), $w$ represents the data-fidelity weight parameter $\lambda$ (\texttt{weight}), and $\alpha = 1.0$ controls the transition stochasticity matching the original DDPM sampling scheme. Guidance is computed using the residual norm (\texttt{guidance='norm'}). Across all tasks, the denoiser $D_{\sigma_t}$ is instantiated using the exact same GS-DRUNet architecture \citep{hurault2022proximal} as in \method and DiffPIR.

\begin{table}[h!]
\centering
\caption{Hyperparameter settings for DPS (\texttt{deepinv} implementation).}
\label{table:hyperparameters_DPS}
\vspace{0.5em}
\begin{tabular}{lccc}
\toprule
Task & Fidelity Weight ($w$) & Steps ($N$) & Noise Schedule \\
\midrule
Inpainting ($\sigma_y = 1/255$) & 4.8 & 300 & \texttt{vp} \\
Inpainting ($\sigma_y = 5/255$) & 4.8 & 300 & \texttt{vp} \\
\bottomrule
\end{tabular}
\end{table}

\subsection{Additional Results}
\subsubsection{Additional Results for MRI}
\label{app:additional_MRI}

First, we recall the Magnetic Resonance Imaging (MRI) formulation. Let $x \in \sC^d$ denote the complex MR image to be reconstructed from undersampled k-space measurements $y \in \sC^d$, obtained via
\begin{equation}
y = \mM\mF x + \rvn,
\end{equation}
where $F \in \sC^{d \times d}$ is the discrete Fourier transform, $\mM \in \{0,1\}^d$ is a binary sampling mask, and $\rvn \sim \mathcal{N}(0, \sigma_y^2 I_d)$ is additive Gaussian noise.
The corresponding data-consistency term is
\begin{equation}
f(x,y) = \frac{1}{2}\|\mM\mF x - y\|_2^2,
\end{equation}
which gradient, required in \method iterations, writes
\begin{equation}
\nabla f(x,y) = \mF^{-1}\mM^\top(\mM\mF x - y).
\end{equation}

Since $F$ is unitary and $\mM$ is a projection, the Lipschitz constantof $\nabla f$ is $L = 1$, so convergence holds for $\tau < 1$.
For \methodP, we need to compute the proximal operator of this data-consistency, which is defined as
\begin{equation}
\operatorname{prox}_{\tau f}(x^k)
=
\arg\min_{x}
\left(
\frac{1}{2}\|\mM\mF x - y\|_2^2
+ \frac{1}{2\tau}\|x - x^k\|_2^2.
\right).
\end{equation}

It admits the closed-form solution
\begin{equation} \operatorname{prox}_{\tau f}(x^k) = \frac{\tau \mF^{-1}y + x^k}{\tau \mM + I}, \end{equation}
where operations are taken component-wise.

\begin{figure}[!ht]
    \centering
    \includegraphics[width=0.84\linewidth]{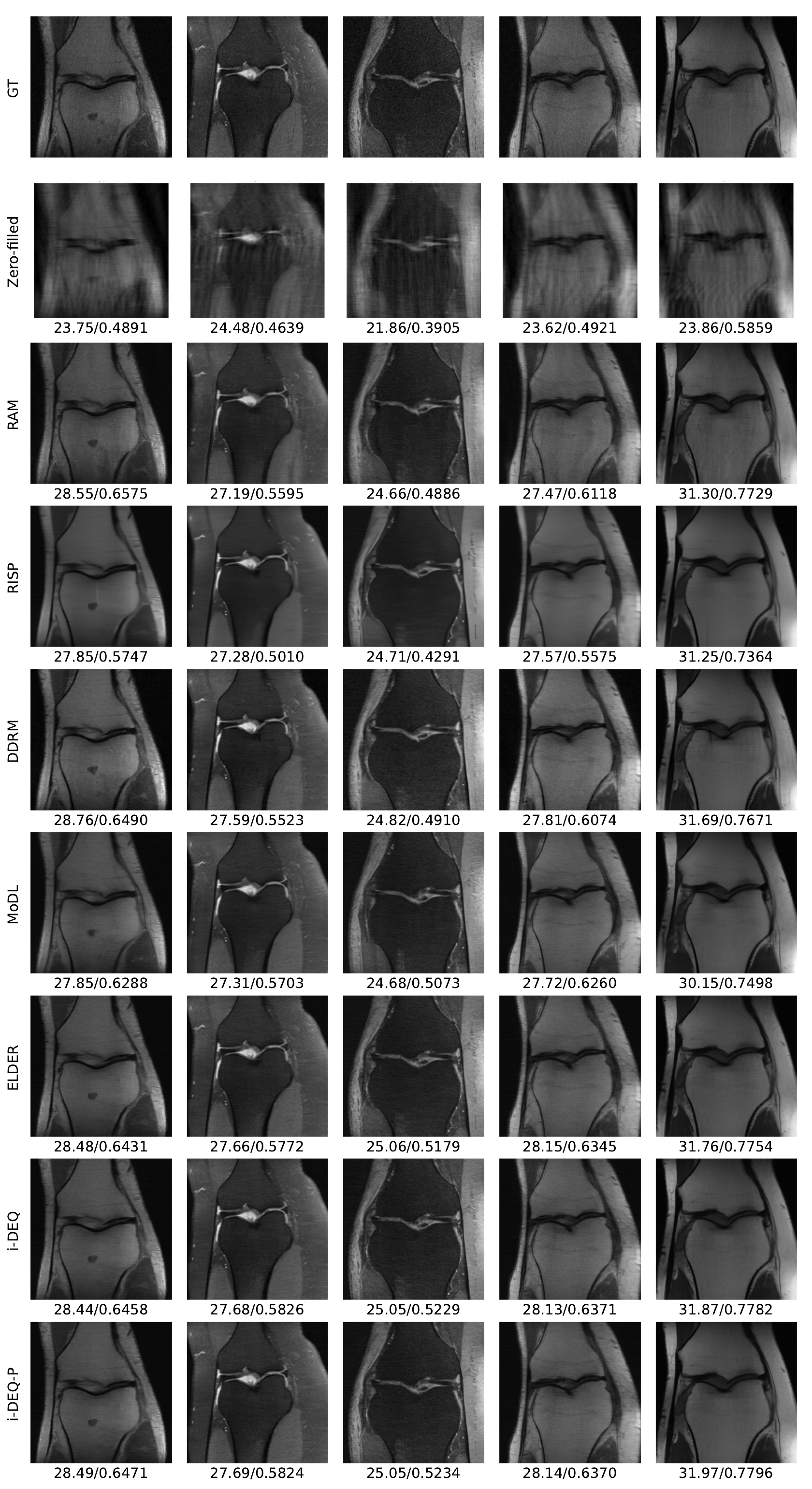}
    \caption{MRI reconstruction results (PSNR/SSIM) for an acceleration factor $\times 8$ and Gaussian noise level $\sigma_y = 1/255$.}
    \label{fig:mri_x8_noise_1}
\end{figure}

Since $x$ is represented in terms of magnitude and phase, with a phase identically equal to zero making $x \in \sR^d$, the network is applied only to the magnitude component of the image while keeping the phase fixed to zero. This choice is convenient because no pretrained denoiser of the GS-Denoiser type is available for complex-valued MRI data. 

\Figref{fig:mri_x8_noise_1} presents reconstruction results for accelerated MRI (acceleration factor $\times 8$) under Gaussian noise with standard deviation $\sigma_y = 1/255$, comparing the evaluated methods.

\subsubsection{Additional Results for Inpainting}
\label{app:additional_inpaint}
Let $x \in \sR^d$ denote the unknown image to be reconstructed from partially observed pixels $y \in \sR^d$. In the inpainting setting, observations are given by
\begin{equation}
y = \mM x + \rvn,
\end{equation}
where $\mM \in \{0,1\}^d$ is a binary masking operator that removes a subset of pixels, and $\rvn \sim \mathcal{N}(0, \sigma_y^2 I_d)$ denotes additive Gaussian noise.
The corresponding data-consistency term reads
\begin{equation}
f(x,y) = \frac{1}{2}\|\mM x - y\|_2^2.
\end{equation}

Its gradient is given by
\begin{equation}
\nabla f(x,y) = \mM^\top(\mM x - y),
\end{equation}
where $\mM$ acts as a diagonal projection operator. Since $\mM$ is an orthogonal projector, $\nabla f$ is 1-Lipschitz, ensuring stability of first-order schemes for $\tau < 1$.
The proximal operator associated with $f$ is defined as
\begin{equation}
\operatorname{prox}_{\tau f}(x^k)
=
\arg\min_{x \in \sR^d}
\left(
\frac{1}{2}\|\mM x - y\|_2^2
+ \frac{1}{2\tau}\|x - x^k\|_2^2
\right).
\end{equation}

It admits a closed-form expression given by
$$\operatorname{prox}_{\tau f}(x^k) = \frac{\tau \mM y + x^k}{\tau \mM + 1},$$ 
where the division is understood component-wise, i.e., it is applied only on observed pixels while leaving missing pixels equal to the current iterate.

\begin{figure}[!ht]
    \centering
    \includegraphics[width=0.84\linewidth]{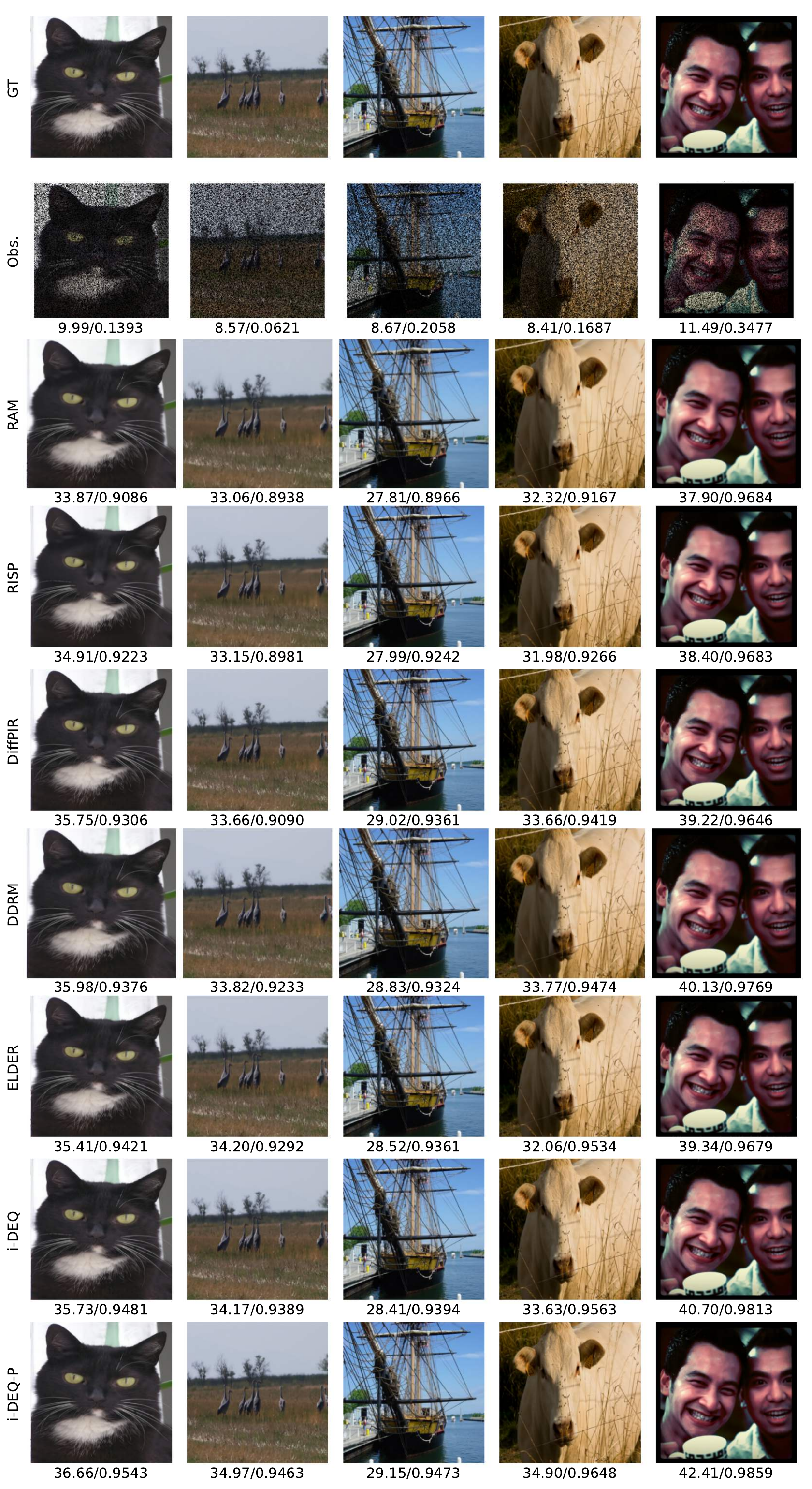}
    \caption{Inpainting results (PSNR/SSIM) with 50\% missing pixels and Gaussian noise level $\sigma_y = 1/255$.}
    \label{fig:inpainting_sigma_1}
\end{figure}

\begin{figure}[!ht]
    \centering
    \includegraphics[width=0.95\linewidth]{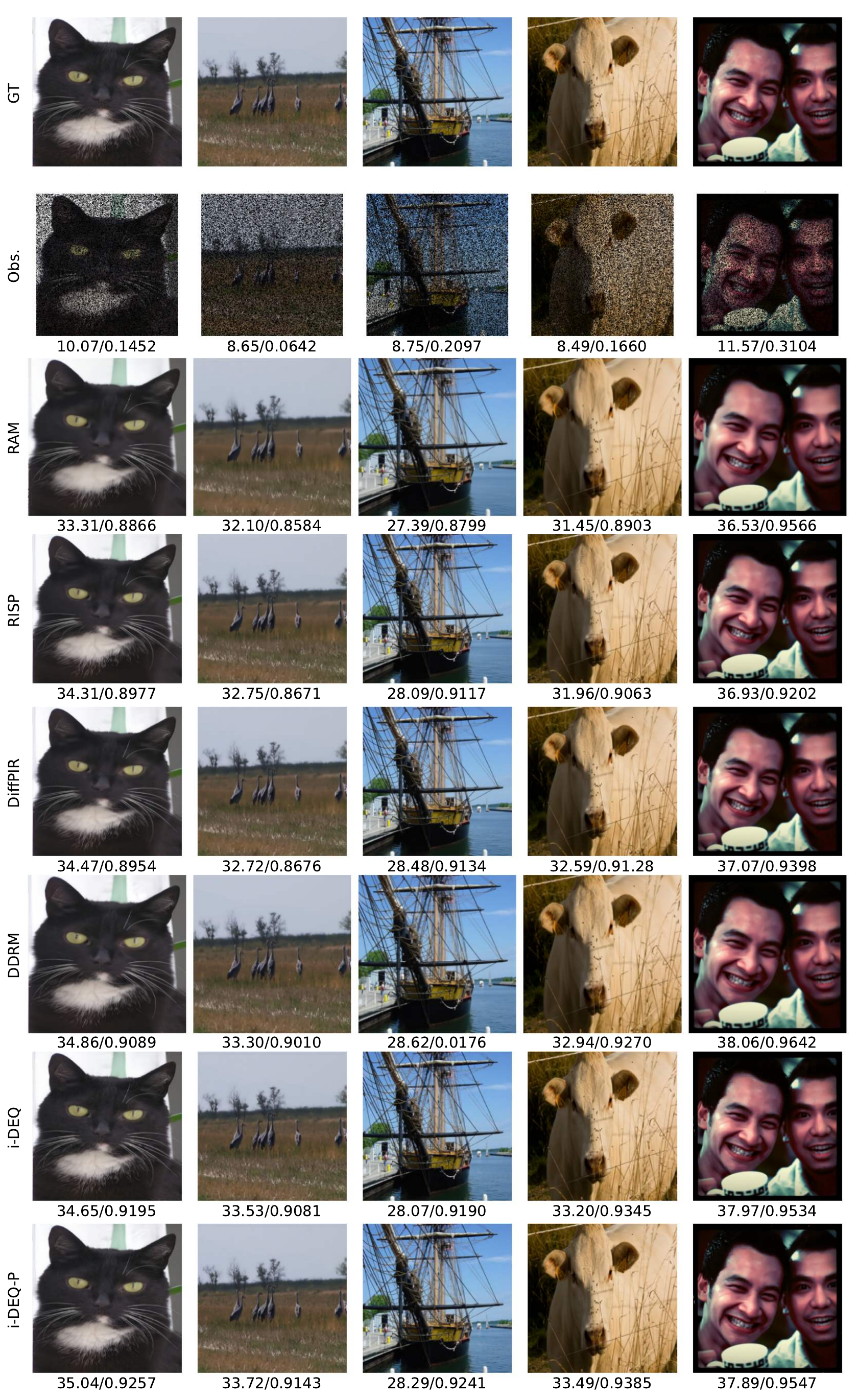}
    \caption{Inpainting results (PSNR/SSIM) with 50\% missing pixels and Gaussian noise level $\sigma_y = 5/255$.}
    \label{fig:inpainting_sigma_5}
\end{figure}

Unlike MRI, the inpainting operator acts directly in the spatial domain, so no transform-domain operations are required. This makes the forward model significantly simpler, while preserving the same optimization structure as in the MRI case.

\twoFigref{fig:inpainting_sigma_1}{fig:inpainting_sigma_5} present additional results for the image inpainting task with 50\% missing pixels under Gaussian noise $\sigma_y = 1/255$ and $\sigma_y = 5/255$.
\subsubsection{Additional Results for Motion Deblurring}
\label{app:additional_deblur}

Let $x \in \sR^d$ denote the clean image to be recovered from a degraded observation $y \in \sR^d$. Under periodic boundary conditions, the motion deblurring model is given by
\begin{equation}
y = \mA x + \rvn,
\end{equation}
where $\mA \in \sR^{d \times d}$ represents a spatial 2D convolution operator with a motion blur kernel $k$, and $\rvn \sim \mathcal{N}(0, \sigma_y^2 I_d)$ denotes additive Gaussian noise. 
Random blur kernels are generated using the Gaussian Process-based trajectory model from \citet{schuler2015learning} via the \texttt{deepinv.physics.generator.MotionBlurGenerator} class with parameters \texttt{psf\_size = 21}, $\sigma = 1.0$, and $l = 0.5$.
The corresponding data-consistency term is defined as
\begin{equation}
f(x,y) = \frac{1}{2}\|\mA x - y\|_2^2.
\end{equation}

By assuming periodic boundary conditions, $\mA$ is diagonalized in the discrete Fourier domain via the 2D Fast Fourier Transform ($\mathcal{F}$). Let $\hat{K} = \mathcal{F}(\mathcal{P}(k))$ denote the complex transfer function of the kernel $k$ after zero-padding and phase-alignment $\mathcal{P}$. The gradient of $f$ is given by
\begin{equation}
\nabla f(x,y) = \mA^\top(\mA x - y) = \operatorname{Re}\left( \mathcal{F}^{-1} \left( |\hat{K}|^2 \odot \mathcal{F}(x) - \overline{\hat{K}} \odot \mathcal{F}(y) \right) \right),
\end{equation}
where $\overline{\hat{K}}$ denotes the complex conjugate of $\hat{K}$, $\odot$ represents the Hadamard product, and $\mA^\top$ corresponds to filtering with the transposed (time-reversed) kernel. Since $\mA$ is a bounded linear operator, $\nabla f$ is $L$-Lipschitz continuous with $L = \|\hat{K}\|_\infty^2$, ensuring numerical stability of first-order optimization schemes for step sizes $\tau < 2/L$.

The proximal operator associated with $f$ is defined as
\begin{equation}
\operatorname{prox}_{\tau f}(x^k)
=
\arg\min_{x \in \sR^d}
\left(
\frac{1}{2}\|\mA x - y\|_2^2
+ \frac{1}{2\tau}\|x - x^k\|_2^2
\right).
\end{equation}

By applying Parseval's identity, this quadratic problem decouples component-wise in the frequency domain, yielding the explicit closed-form expression
\begin{equation}
\operatorname{prox}_{\tau f}(x^k) = \operatorname{Re}\left( \mathcal{F}^{-1} \left( \frac{\mathcal{F}(x^k) + \tau \overline{\hat{K}} \odot \mathcal{F}(y)}{1 + \tau |\hat{K}|^2} \right) \right),
\end{equation}
where the division is applied element-wise across frequency coordinates.

Unlike inpainting, the motion deblurring operator operates globally in the spatial domain via a shift-invariant blur kernel, requiring Fourier-domain computations to efficiently perform convolutions and inversions while maintaining the same optimization pipeline.

\begin{figure}[!ht]
    \centering
    \includegraphics[width=\linewidth]{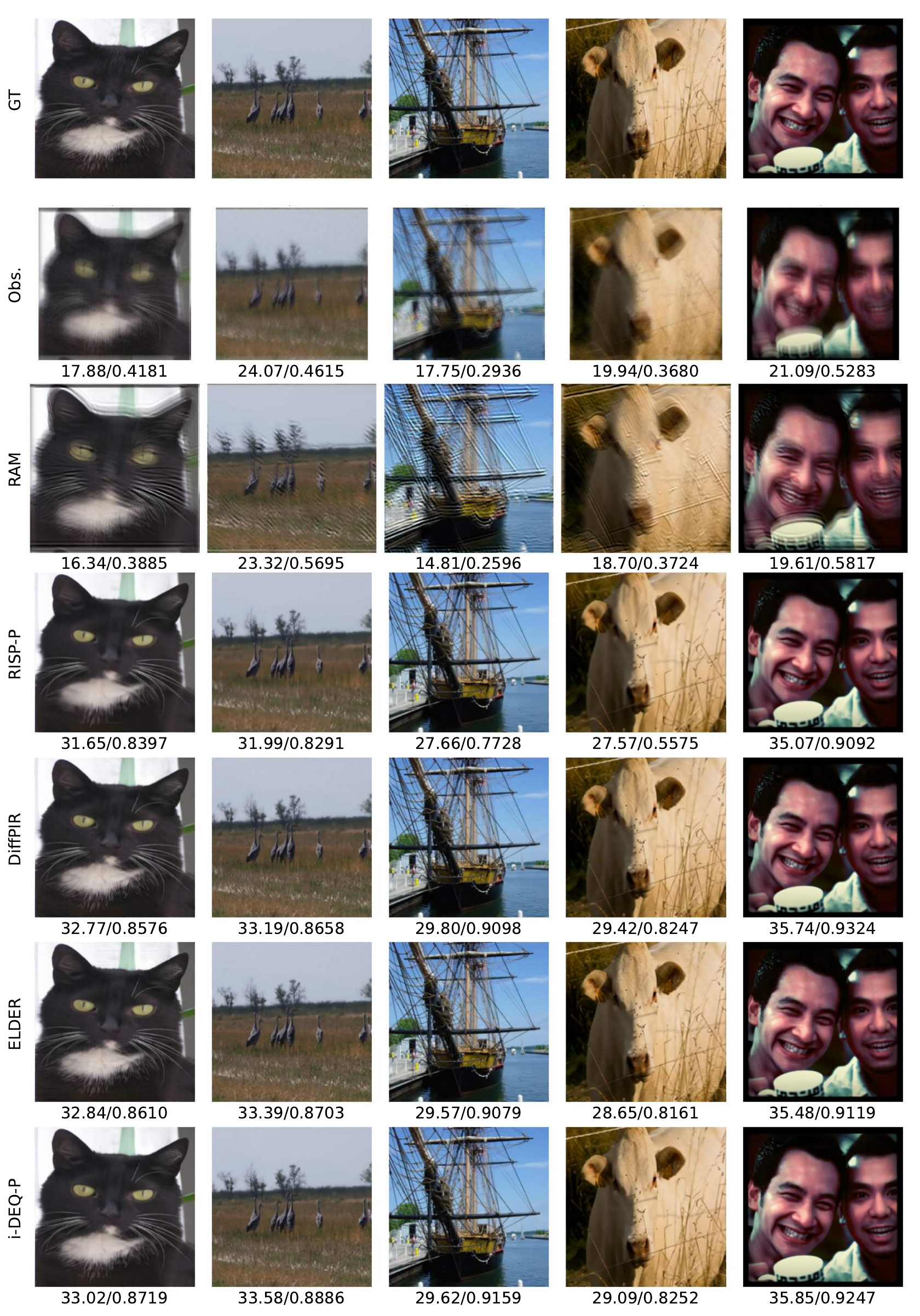}
    \caption{Motion debluring results (PSNR/SSIM) at noise level $\sigma_y = 5/255$.}
    \label{fig:deblur_sigma_5}
\end{figure}

\Figref{fig:deblur_sigma_5} presents additional results for the motion deblurring task.

\subsubsection{Additional Results for Rician Denoising}
\label{app:additional_rician}

Let $n_{\text{real}}$ and $n_{\text{imag}}$ be two independent Gaussian random variables with zero mean and standard deviation $\sigma_y$. The observed magnitude data $y$ is given by
\begin{equation}
y = \sqrt{(x + n_{\text{real}})^2 + n_{\text{imag}}^2}.
\end{equation}

The conditional distribution of $y$ given $x$ follows a Rician distribution, whose probability density function is
\begin{equation}
p(y \mid x) = \frac{y}{\sigma_y^2} \exp\!\left(-\frac{x^2 + y^2}{2\sigma_y^2}\right)
I_0\!\left(\frac{x y}{\sigma_y^2}\right),
\end{equation}
where $I_0$ (resp $I_1$) denotes the modified Bessel function of the first kind of order zero (resp one).

From this model, the associated data-consistency term (up to additive constants) can be written as
\begin{equation}
f(x,y) = \frac{x^2}{2\sigma_y^2} - \log I_0\!\left(\frac{x y}{\sigma_y^2}\right).
\end{equation}

Its derivative is given by
\begin{equation}
\nabla f(x,y) = \frac{x}{\sigma_y^2} - \frac{y}{\sigma_y^2} B\!\left(\frac{x y}{\sigma_y^2}\right),
\end{equation}
where $B(t) = \frac{I_1(t)}{I_0(t)}$.\\

However,  the associated proximal operator required to apply \methodP does not admit a closed-form expression due to the presence of the modified Bessel function. 
It can be approximated using an iterative reweighted $\ell_1$ minimization scheme (Iteratively Reweighted Least Absolute Deviations / $\ell_1$-type method), which can be interpreted as a majorization-minimization strategy for the non-convex Rician log-likelihood involving modified Bessel functions~\citep{daubechies2003iterativethresholdingalgorithmlinear,KOAY2006317}.
Despite its nonlinearity, the function $f$ is smooth and has Lipschitz continuous gradient and Hessian~\citep[Lemma 12]{RISP}, which enables the use of RISP and the proposed \method.

\begin{figure}[!ht]
    \centering
    \includegraphics[width=.95\linewidth]{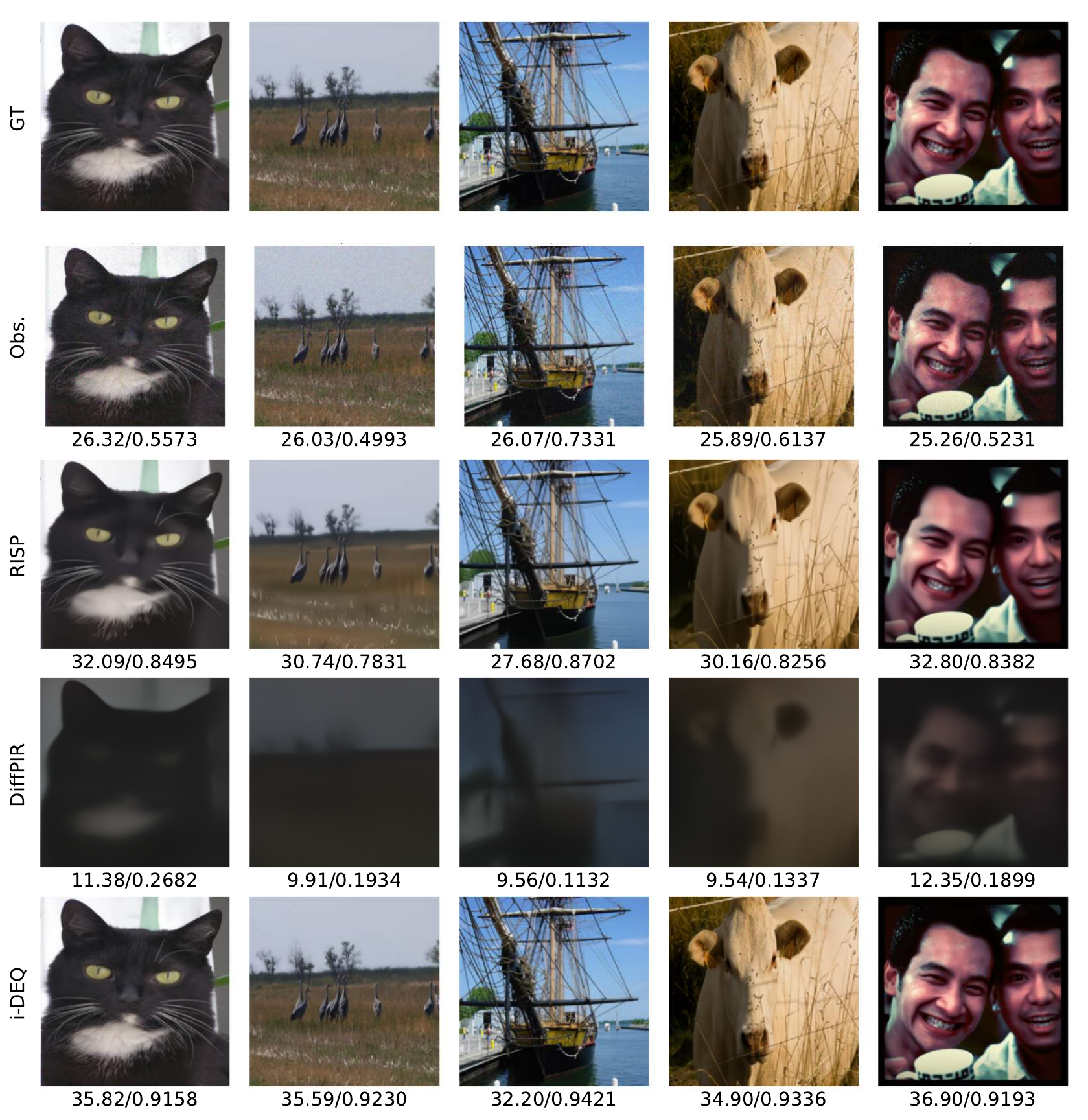}
    \caption{Rician denoising results (PSNR/SSIM) at noise level $\sigma_y = 0.05$.}
    \label{fig:rician_sigma_1}
\end{figure}

\begin{figure}[!ht]
    \centering
    \includegraphics[width=\linewidth]{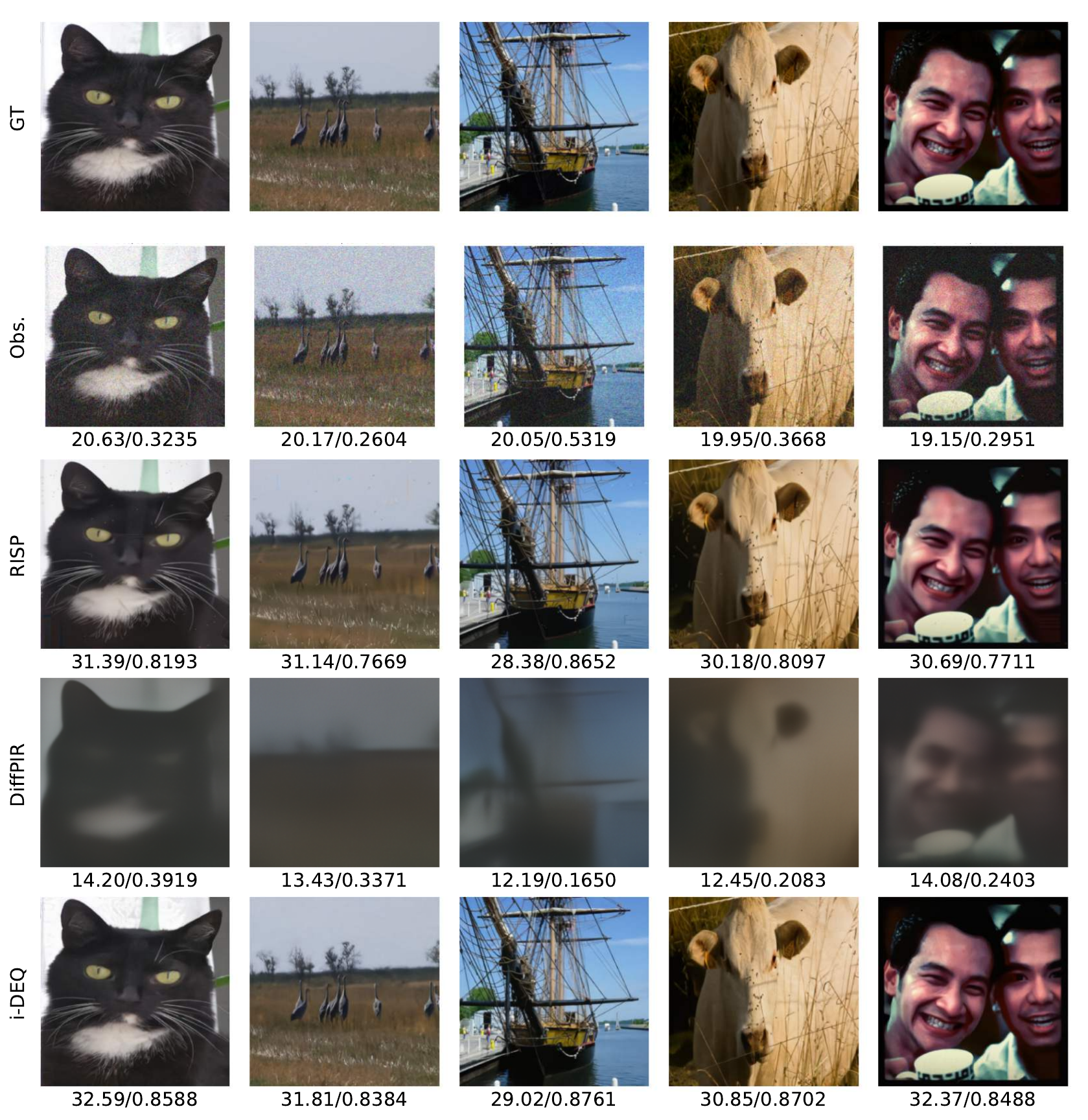}
    \caption{Rician denoising results (PSNR/SSIM) at higher noise level $\sigma_y = 0.1$.}
    \label{fig:rician_sigma_5}
\end{figure}

\twoFigref{fig:rician_sigma_1}{fig:rician_sigma_5} present additional results for Rician denoising at noise levels $\sigma_y = 0.05$ and $\sigma_y = 0.1$.

\FloatBarrier

\newpage

\section*{NeurIPS Paper Checklist}
\begin{enumerate}

\item {\bf Claims}
\item[] Question: Do the main claims made in the abstract and introduction accurately reflect the paper's contributions and scope?
\item[] Answer: \answerYes{}
\item[] Justification: The claims stated in the abstract and introduction accurately reflect the scope and contributions of the paper. The theoretical and experimental results are consistent with the presented methodology, and the reported performance improvements are supported by the conducted experiments without overstatement beyond the evaluated settings.

\item {\bf Limitations}
\item[] Question: Does the paper discuss the limitations of the work performed by the authors?
\item[] Answer: \answerYes{}
\item[] Justification:
The paper explicitly states limitations regarding (i) computational efficiency, since i-DEQ remains slower than diffusion methods such as DiffPIR, and (ii) robustness, as training can still become unstable in high-noise regimes. These points directly reflect acknowledged constraints on scalability and stability.

\item {\bf Theory assumptions and proofs}
    \item[] Question: For each theoretical result, does the paper provide the full set of assumptions and a complete (and correct) proof?
    \item[] Answer: \answerYes{}
    \item[] Justification: All assumptions are clearly stated and discussed in the paper. The main proposition is reused from a cited work, where the complete proof is provided in detail.

    \item {\bf Experimental result reproducibility}
    \item[] Question: Does the paper fully disclose all the information needed to reproduce the main experimental results of the paper to the extent that it affects the main claims and/or conclusions of the paper (regardless of whether the code and data are provided or not)?
    \item[] Answer: \answerYes{}
    \item[] Justification: Code and data are provided. Hyperparameters are detailed in the appendix, and fixed random seeds ensure reproducibility. Minor differences may arise due to variations in GPU hardware.

    \item {\bf Open access to data and code}
    \item[] Question: Does the paper provide open access to the data and code, with sufficient instructions to faithfully reproduce the main experimental results, as described in supplemental material?
    \item[] Answer: \answerYes{}
    \item[] Justification: Code and data will be made publicly available upon publication and are provided with the submission, with sufficient documentation to reproduce the main results. Hyperparameters and implementation details are provided, and fixed random seeds ensure reproducibility. Minor discrepancies may occur due to hardware differences.

    \item {\bf Experimental setting/details}
    \item[] Question: Does the paper specify all the training and test details (e.g., data splits, hyperparameters, how they were chosen, type of optimizer) necessary to understand the results?
    \item[] Answer: \answerYes{}
    \item[] Justification: The paper specifies the experimental setup, including data splits, hyperparameters, and optimization details. Additional information is provided in the appendix and supplemental material, ensuring that the results can be properly understood and reproduced.

\item {\bf Experiment statistical significance}
\item[] Question: Does the paper report error bars suitably and correctly defined or other appropriate information about the statistical significance of the experiments?
\item[] Answer: \answerNo{}
\item[] Justification: The paper does not provide detailed error bars or statistical significance analysis. However, the test and validation sets are relatively small, making standard deviations less informative, while the mean PSNR is sufficient to assess overall performance. Moreover, there is high intrinsic variability across images, and all methods are evaluated on the same data, ensuring a fair comparison under shared sources of variability. In addition, a sufficiently large number of reconstructed examples is provided in the appendix to illustrate the range of qualitative outcomes and to support the observed variability in results.

\item {\bf Experiments compute resources}
\item[] Question: For each experiment, does the paper provide sufficient information on the computer resources (type of compute workers, memory, time of execution) needed to reproduce the experiments?
\item[] Answer: \answerYes{}
\item[] Justification: The paper specifies the compute resources used, including the type of hardware (GPU), and provides sufficient details to reproduce the experiments. Execution time are detailed.
    
\item {\bf Code of ethics}
\item[] Question: Does the research conducted in the paper conform, in every respect, with the NeurIPS Code of Ethics \url{https://neurips.cc/public/EthicsGuidelines}?
\item[] Answer: \answerYes{}
\item[] Justification: The research complies with the NeurIPS Code of Ethics. It does not involve human subjects or sensitive data, and relies on publicly available datasets with appropriate licenses. No personally identifiable information is used. The proposed methods are evaluated in a controlled setting and do not raise concerns related to safety, security, or misuse. Potential biases are inherent to the datasets and affect all compared methods equally. The work does not involve applications with direct societal or environmental risks.

\item {\bf Broader impacts}
\item[] Question: Does the paper discuss both potential positive societal impacts and negative societal impacts of the work performed?
\item[] Answer: \answerYes{}
\item[] Justification: The paper discusses both positive and negative societal impacts in a dedicated section.
    
\item {\bf Safeguards}
\item[] Question: Does the paper describe safeguards that have been put in place for responsible release of data or models that have a high risk for misuse (e.g., pre-trained language models, image generators, or scraped datasets)?
\item[] Answer: \answerNA{}
\item[] Justification: The proposed work does not involve high-risk or dual-use models such as generative models or systems with direct misuse potential. It focuses on image reconstruction from corrupted measurements and does not include deployment-oriented components. Consequently, no specific safeguards for model or data release are required.

\item {\bf Licenses for existing assets}
\item[] Question: Are the creators or original owners of assets (e.g., code, data, models), used in the paper, properly credited and are the license and terms of use explicitly mentioned and properly respected?
\item[] Answer: \answerYes{}
\item[] Justification: The paper properly credits the original sources of all datasets used and cites the corresponding publications as required. Both datasets are standard benchmarks in the literature and are used in compliance with their respective licenses and terms of use. The specific dataset versions are clearly identified in the experimental section, and the usage follows the conventions established in prior publications.

\item {\bf New assets}
\item[] Question: Are new assets introduced in the paper well documented and is the documentation provided alongside the assets?
\item[] Answer: \answerYes{}{}
\item[] Justification: The paper releases code along with pre-trained model weights to ensure reproducibility and facilitate evaluation of the proposed method. The provided repository includes documentation describing how to install dependencies, run training and inference procedures, and reproduce the reported results.

\item {\bf Crowdsourcing and research with human subjects}
\item[] Question: For crowdsourcing experiments and research with human subjects, does the paper include the full text of instructions given to participants and screenshots, if applicable, as well as details about compensation (if any)? 
\item[] Answer: \answerNA{}
\item[] Justification: The paper does not involve crowdsourcing experiments or research with human subjects. All experiments are conducted on publicly available image datasets using computational evaluation only.

\item {\bf Institutional review board (IRB) approvals or equivalent for research with human subjects}
\item[] Question: Does the paper describe potential risks incurred by study participants, whether such risks were disclosed to the subjects, and whether Institutional Review Board (IRB) approvals (or an equivalent approval/review based on the requirements of your country or institution) were obtained?
\item[] Answer: \answerNA{}
\item[] Justification: The paper does not involve human subjects or any form of data collection from participants. Therefore, no IRB approval or equivalent ethical review was required.

\item {\bf Declaration of LLM usage}
\item[] Question: Does the paper describe the usage of LLMs if it is an important, original, or non-standard component of the core methods in this research? Note that if the LLM is used only for writing, editing, or formatting purposes and does \emph{not} impact the core methodology, scientific rigor, or originality of the research, declaration is not required.
\item[] Answer: \answerNA{}
\item[] Justification: The core methodology does not involve the use of large language models as part of the proposed approach. Any potential LLM usage is limited to auxiliary purposes such as writing or editing and does not affect the scientific content or experimental results.

\end{enumerate}
\end{document}